\documentclass[11pt]{article}
\usepackage{amssymb}
\usepackage{graphicx}

 {\catcode `\@=11 \global\let\AddToReset=\@addtoreset}
\AddToReset{equation}{section}

\newtheorem{theorem}{Theorem}[section]

\newtheorem{corollary}{\bf Corollary}[section]

\newtheorem{@definition}{\sc Definition}[section]

\newtheorem{@remark}{\sc Remark}[section]

\newtheorem{@example}{\sc Example}[section]

\addtocounter{theorem}{2}

\def\mathsf{\bf}

\def\d{\mathrm d}

\def\E{\mathrm E}

\def\text{\mbox}

\def\g{\gamma}
\def\a{\alpha}
\def\b{\beta}

\def\text{\mbox}

\newcommand{\beq}{\begin{equation}}
\newcommand{\eeq}{\end{equation}}
\newcommand\beqn{\begin{displaymath}}  
\newcommand\eeqn{\end{displaymath}}

\newcommand{\halmos}{\vspace{3mm} \hfill \mbox{$\Box$}\\[2mm]}

\begin{document}

%

\title{Several new tail index estimators
\footnotemark[0]\footnotetext[0]{  \scriptsize
The research was supported by  Research Council of Lithuania, grant
No. \ MIP-076/2013} }

\author{Vygantas Paulauskas$^{\text{\small a}}$, Marijus Vai{\v c}iulis$^{\text{\small b,a}}$ \\
{\tiny $^{\text{a}}$ Department of Mathematics and Informatics,Vilnius University, Naugarduko st. 24, LT-03225 Vilnius, Lithuania}\\
{\tiny $^{\text{b}}$  Institute of Mathematics and Informatics, Vilnius University, Akademijos st. 4, LT-08663 Vilnius, Lithuania}\\
 }


\maketitle

\begin{abstract}
In the paper we propose some new class of functions which is used to
construct tail index estimators.  Functions from this new class is non-monotone in general, but presents a product of two monotone  functions: the power function and the logarithmic function, which plays
essential role in the classical Hill estimator. Introduced   new
estimators have better  asymptotic performance comparing with  the
Hill estimator and other  popular estimators over all range of the
parameters present in the second order regular variation condition.
Asymptotic normality of the  introduced estimators is proved, and
comparison (using asymptotic mean square error) with other
estimators of the tail index is provided. Some preliminary simulation results are presented.
\end{abstract}

\vfill
\eject

\section{Introduction and formulation of results }

From the first papers of Hill and Pickands (see \cite{Hill} , \cite{Pickands}), devoted to the estimation of the tail index (or, more generally, extreme value index), most of statistics constructed for this aim were based on order statistics and logarithmic function. Suppose we have a sample $X_1, X_2, \dots , X_n$, considered as independent, identically distributed (i.i.d.) random variables with a distribution function
d.f. $F$ satisfying the following relation for large $x$:
\begin{equation}\label{cond1}
{\bar F}(x):=1-F(x)=x^{-1/\g}L(x).
\end{equation}
Parameter $\g$ is called extreme value index (EVI) and  $\a:=1/\g >0$ is called the tail index,   $L(x) >0$, for all
$x>0$, and $L$ is a slowly varying at infinity function:
$$
\lim_{x\to \infty}{L(tx) \over L(x)}=1.
$$
In the paper we consider the case $\g >0$ only.
Denote
$U(t)=F^{\leftarrow}\left(1-(1/t)\right)^{\leftarrow}$, $t \ge 1$,
 where $W^{\leftarrow}:I \to
\mathbb{R}$ is the left continuous inverse function of a monotone
function $W$, defined by the relation
$W^{\leftarrow}(t)=\sup\left\{x: \ W(x) \le t\right\}$, $t \in I$. It is well-known that in the case $\g>0$  assumption (\ref{cond1}) is equivalent to the following one: for all $x>0$,
\begin{equation}\label{cond2}
\lim_{t\to \infty} \frac{U(tx)}{U(t)}=x^{\g},
\end{equation}
i.e., the quantile function $U(t)$ varies regularly with index $\g$.
Let $X_{n,1}\le X_{n,2}\le \cdots \le
X_{n,n}$ denote the order statistic of $X_1, \dots ,X_n.$  Taking some part of the largest values from the sample and the logarithmic function a statistician can form various statistics. In this way one can get Hill and Pickands estimators, moment and moment ratio estimators which are well-known  and  deeply investigated. The heuristic behind this approach  (the threshold-over-peaks (TOP) phenomenon and maximum likelihood) is also given in many papers and monographs, therefore we do not provide it here.

There were estimators, based on a little bit different idea: the sample is divided into blocks and in each block the ratio of  two largest values is taken. Then the linear function  $f(x)=x$ instead of logarithmic one is applied to these ratios. Estimators, based on this idea were constructed in  \cite{Paul8} and  \cite{Paul17}.
The next step was to include the linear and logarithmic functions into some parametric family of functions, and, considering estimators based on block scheme, this was done in \cite{Paul18}, taking the family of functions, defined for $x\ge 1$
\begin{equation} \label{f_r}
f_r(x)= \left\{
  \begin{array}{ll}
    \frac{1}{r}\left(x^r -1 \right), & \hbox{$r<\alpha$, $r \not=0$,} \\
   \ln x, & \hbox{$r=0$.}
  \end{array}
\right.
\end{equation}
In \cite{Paul16} this  family of functions  was applied for order statistic. This was done by introducing the statistics
$$
   H_n^{(\ell)}(k,r)=\frac{1}{k} \sum_{i=0}^{k-1} f_r^{\ell} \left( \frac{X_{n-i,n}}{X_{n-k,n}}\right), \quad \ell=1,2,
$$
and some combinations, formed from these statistics.
Here  $k$ is some number satisfying $1\le k < n$, and in the EVI
estimation $k$ is chosen as a function of $n$ (thus, strictly speaking
we should denote it by $k(n)$). In this way the generalizations of the Hill, the moment, and the  moment ratio estimators (we shall write the expressions of these estimators a little bit later; to denote these estimators we shall use the letters $H, M, MR$ and for generalizations we add the letter $G$) were obtained, for example, the generalized Hill estimator (GHE) is defined as
\begin{equation}\label{ghe}
\hat{\gamma}_n^{(GH)} (k,r)=\frac{H_n(k,r)}{1+r \cdot H_n(k,r)}
\end{equation}
and the Hill  estimator is obtained by taking $r=0, \ \hat{\gamma}_n^{(H)} (k):=\hat{\gamma}_n^{(GH)} (k,0)$.

The main goal of this paper is to introduce another parametric family of functions, which has the same property that includes logarithmic function, and to construct new estimators using this family.
For $x\ge 1$ let us consider functions
$$
g_{r,u}(x)=x^r \ln^u(x),
$$
where parameters $r$ and $u$ can be arbitrary real numbers, but for our purposes, connected with consistency, we shall require $\g r<1$ and $u>-1$.
Moreover, mainly we shall consider only integer values of parameter $u$. The family $\{g_{r, u}\}$, similarly to $\{f_r\}$, contains logarithmic function (with $r=0$), but, contrary to $\{f_r\}$, contains logarithmic function for any value of parameter $r$ (if $u \ne 0$). Also let us note that for $r\ge 0$ function $g_{r,u}$ is monotone for all values of $u$, while for $r<0$ and $u>0$ there is no monotonicity.

Using these functions we can form statistics, similar to $ H_n^{(\ell)}(k,r)$:
$$
G_n(k,r,u)=\frac{1}{k} \sum_{i=0}^{k-1} g_{r,u} \left(
\frac{X_{n-i,n}}{X_{n-k,n}}\right).
$$
The above mentioned Hill estimator, moment estimator (introduced in \cite{Dekkers}) and  the moment ratio estimator (introduced in \cite{Danielsson})
can be expressed via statistics $G_n(k,0,u)$, $u=1,2$ as follows
\begin{eqnarray*}
\hat{\gamma}^{H}_n(k)&=& G_n(k,0,1), \\
\hat{\gamma}^{M}_n(k)&=& G_n(k,0,1)+\frac{1}{2} \left\{1-\left( \frac{G_n(k,0,2)}{G_n^2(k,0,1)}-1\right)^{-1} \right\}, \\
\hat{\gamma}^{MR}_n(k)&=& \frac{G_n(k,0,2)}{2 G_n(k,0,1)}.
\end{eqnarray*}
Let us note, that, due to the expressions of functions $f_r$ and $g_{r,u}$, we can express statistic $H_n(k,r)$ via statistics $G_n(k,r,u)$:
\begin{equation} \label{sarysis}
H_n(k,r)= \left\{
  \begin{array}{ll}
    \left( G_n(k,r,0)-1\right)/r,   &  \hbox{ $r \not=0$,} \\
   G_n(k,0,1), & \hbox{$r=0$,}
  \end{array}
\right.
\end{equation}
and there is continuity with respect to $r$ in this relation, since it is easy to see that $\lim_{r\to 0}\left( G_n(k,r,0)-1\right)/r = G_n(k,0,1).$
 Taking into account that $H_n^{(2)}(k,r)$ can be expressed via $H_n(k,r)$ (see (3.2) in \cite{Paul16}), all estimators, which were introduced in \cite{Paul16}, can be written by means of statistics $G_n(k,r,u)$ only.

One more estimator of the parameter $\g$, which can be written as a function of statistics $ H_n^{(\ell)}(k,r)$, therefore also by statistics $G_n(k,r,u),$
very recently was introduced in \cite{Beran} (see also \cite{Henry} and \cite{Stehlik}).
It is named as the harmonic moment tail index estimator (HME, this abbreviation is used in \cite{Beran}),
and for $\beta >0, \ \b \ne 1,$ it is defined by formula
\begin{equation}\label{defHME}
\g_{n}^{HM}(k, \b) = \frac{1}{\beta-1} \left( \left[ k^{-1} \sum_{i=0}^{k-1} \left( \frac{X_{n-k,n}}{X_{n-i,n}}\right)^{\beta-1} \right]^{-1}-1\right),
\end{equation}
while for $\b =1$ it is defined as a limit as $\b \to 1$. It is easy to see that
$$
\g_{n}^{HM}(k, 1) = k^{-1} \sum_{i=0}^{k-1} \ln \left(\frac{X_{n-k,n}}{X_{n-i,n}}\right),
$$
i.e., we get the Hill estimator. But very simple transformation shows that, denoting $\b=1-r$, we have
 $$
  \g_{n}^{HM}(k, 1-r)=\frac{1- (r \cdot H_n^{(1)} (k,r) +1)^{-1}}{r} = \frac{ H_n^{(1)} (k,r)}{ r \cdot H_n^{(1)} (k,r) +1 }= \gamma_n^{GH} (k,r).
 $$
 This means, that the HME is exactly the GHE, only additional (tuning) parameters $\b$ (for HME) and $r$ (for GHE) have different ranges.
It is strange why there is requirement $\b >0$ in the definition of HME, which corresponds to the condition $r<1$ for GHE. Of course, both statistics can be
defined for all real  $\b$ and $r$, but in order to have consistent estimators of the  parameter $\g$ we must have some restrictions on these tuning parameters.
 In \cite{Paul16} the consistency of GHE is proved under condition $\g r<1$ (or $r<\a$), which is equivalent to condition $\b>1-\g^{-1}$, while in \cite{Beran}
 the consistency of HME is stated (there is only the sketch of the proof) without any additional condition on $\b$, so one can think that  HME is consistent for
 all $\b>0$. In \cite{Paul16} it was explained why the GHE is not consistent if $r>\a$, and from the expression (\ref{defHME}) intuitively it is clear that for
 $\b<1-\g^{-1}$ we have $\g_{n}^{HM}(k, \b){\buildrel \rm P \over \longrightarrow}(1-\b)^{-1}$. This discussion even shows that the requirement $\b>0$ is not natural: for $0<\g <1$ we have $1-\g^{-1}<0$, thus for
$1-\g^{-1}<\b <0$ we still have consistency of HME, while if $\g>1$, then in the interval $0<\b <1-\g^{-1}$ there is no consistency. The asymptotic normality of HME in \cite{Beran} is proved under the requirement $\b>1-(2\g)^{-1}$, which is equivalent to the requirement $2\g r<1$, used in \cite{Paul16}.

Our first result shows what quantities are estimated by the introduced statistics $G_n(k,r,u)$.

\begin{theorem}\label{thm1}
Suppose that $X_1,\dots, X_n$ are i.i.d. nonnegative random variables whose quantile function $U$ satisfies
condition (\ref{cond2}). Let $\g r<1$ and $u>-1$. Let us suppose that  a sequence $k=k(n)$ satisfies conditions
\begin{equation}\label{cond3}
    k(n) \to \infty, \quad n/k(n)\to \infty, \quad {\rm as} \ n \to \infty.
\end{equation}
Then we have
\begin{equation}\label{res1}
G_n(k,r,u) {\buildrel \rm P \over \to}\xi(r, u):= \frac{\gamma^{u} \Gamma(1+u)}{(1-\gamma r)^{1+u}},  \quad {\rm as} \ n \to \infty.
\end{equation}
Here ${\buildrel \rm P \over \to}$ stands for the convergence in probability and $\Gamma(u)$ denotes Euler's gamma function.
\end{theorem}
The following corollary allows to proof the consistency of estimators, expressed as a function of statistics $G_n(k,r,u)$ with different $r$ and $u.$
\begin{corollary} \label{cor1}
Suppose that $X_1,\dots, X_n$ are i.i.d. nonnegative random variables whose quantile function $U$ satisfies
condition (\ref{cond2}). Let $\gamma r_{j}<1$ and  $u_{j}>-1$, $j=1,2,\dots,\ell$.  Let us suppose that  a sequence $k=k(n)$ satisfies (\ref{cond3}).
Let $\chi(t_1, \dots, t_{\ell}): (0,\infty)^{\ell} \to (0,\infty)$ be a continuous function. Then
\begin{equation} \label{cor1-res}
\chi\left( G_n(k,r_1,u_1), \dots, G_n(k,r_{\ell},u_{\ell}) \right) {\buildrel \rm P \over \to}
\chi\left(
\xi(r_1, u_1), \dots, \xi(r_{\ell}, u_{\ell})\right),
\end{equation}
as $n \to \infty$.
\end{corollary}

The relation (\ref{cor1-res})  gives us many possibilities to form consistent estimators of $\g$ using
statistics $\ G_n(k,r,u)$ with different $r,u$.
  Since it is not clear which combinations are good ones, we decided to restrict ourselves with two most simple statistics $G_n(k,r, 0), \ G_n(k,r, 1)$ (that is, $u=0$  and $u=1$)and to consider the following three estimators of $\g>0$
\begin{eqnarray}
  \hat{\gamma}^{(1)}_n(k,r) &=&
  \left\{
    \begin{array}{ll}
      (G_n(k,r,0)-1)(r G_n(k,r,0))^{-1}, & \hbox{$r\not = 0$,} \\
      G_n(k,0,1), & \hbox{$r=0$,}
    \end{array}
  \right.
  \nonumber \\
\hat{\gamma}^{(2)}_n(k,r) &=& \frac{2 G_n(k,r,1)}{2 r G_n(k,r,1) +1+ \sqrt{4r
  G_n(k,r,1)+1}}, \label{statgrnew2} \\
\hat{\gamma}^{(3)}_n(k,r)&=&
\left\{
  \begin{array}{ll}
    (r G_n(k,r,1) -G_n(k,r,0)+1)(r^2 G_n(k,r,1) )^{-1}, & \hbox{$r \not = 0$,} \\
    \hat{\gamma}^{MR}_n(k), & \hbox{$r=0$.}
  \end{array}
\right. \nonumber
\end{eqnarray}
One can note, that the estimator $\hat{\gamma}^{(1)}_n(k,r)$ is exactly the GHE, given in (\ref{ghe}), only expressed via statistics  $G_n(k,r,u)$. For us it will be convenient to use this notation for GHE, since we shall compare these two new constructed estimators not with the Hill  estimator (which earlier was like a "benchmark" for other estimators), but with the GHE $\hat{\gamma}^{(1)}_n(k,r)$. Since $\hat{\gamma}^{(2)}_n(k,0)=\hat{\gamma}^{(1)}_n(k,0)$, the second estimator presents another generalization of the Hill estimator, while the third estimator gives us the generalized moment ratio estimator.

The main step in proving asymptotic normality of the introduced estimators $\hat{\gamma}^{(j)}_n(k,r), j=1,2,3,$ is to prove asymptotic normality for statistics $G_n(k,r, 0), \ G_n(k,r, 1)$. As usual, in order to get asymptotic normality for estimators the so-called second-order regular variation (SORV) condition, in one or another form, is assumed. In this paper we shall use the SORV condition formulated by means of the function $U$. We  assume that there exists a measurable
function $A(t)$ with the constant sign near infinity, which is not identically zero, and $A(t) \to 0$ as $t \to \infty$,
such that
\begin{equation}\label{cond4}
    \lim_{t \to \infty} \frac{\frac{U(tx)}{U(t)}-x^{\gamma}}{A(t)}=x^{\gamma} \frac{x^{\rho}-1}{\rho}
\end{equation}
for all $x > 0$. Here $\rho < 0$ is the so-called second order parameter. It is known that (\ref{cond4})
implies that the function $|A(t)|$ varies regularly with index $\rho$.

Let us denote $d_r(k)=1-k\gamma r$.

\begin{theorem} \label{thm2}
Suppose that $X_1,\dots, X_n$ are i.i.d. nonnegative random variables whose quantile function $U$ satisfies
condition (\ref{cond4}). Suppose that  $\gamma r<1/2$   and   that  the sequence $k=k(n)$ satisfies
(\ref{cond3}) and
\begin{equation}\label{cond5}
    \lim_{n\to \infty} \sqrt{k} A\left( \frac{n}{k}\right)=\mu \in (-\infty, +\infty).
\end{equation}
Then, as $n \to \infty$,
\begin{equation}\label{res2}
    \sqrt{k} \left(G_n(k,r,0)-\xi(r,0), G_n(k,r,1)-\xi(r,1) \right) {\buildrel \rm d \over \to} \mu \left(\nu^{(1)}(r), \nu^{(2)}(r)\right)+
    \left(W^{(1)}, W^{(2)}\right),
\end{equation}
 where ${\buildrel \rm d \over \to}$ stands for the convergence in distribution and quantities $\nu^{(j)}(r), \ j=1,2$ are as follows
\begin{equation} \label{nu}
\nu^{(1)}(r) =\frac{r}{d_r(1)(d_r(1)-\rho)}, \quad \nu^{(2)}(r)
=\frac{1-\rho-\gamma^2 r^2}{(d_r(1))^2(d_r(1)-\rho)^2}.
\end{equation}
Here $\left(W^{(1)},W^{(2)}\right)$ is zero mean Gaussian random vector with variances
$\E\left(W^{(j)}\right)^2=s_j^2 (r)$, $j=1,2$ and covariance $\E\left(W^{(1)}W^{(2)} \right)=s_{12}(r)$, where
\begin{equation} \label{s}
s_1^2(r)=\frac{\gamma^2 r^2}{d_r(2)d_r^2(1)}, \quad s_2^2 (r) = \frac{
\gamma^2(d_r(2)+2\gamma^4 r^4)} {d_r^3(2)d_r^4(1)}, \quad
s_{12}(r)=\frac{\gamma^2 r (d_r(1)-\gamma^2 r^2)}{d_r^2(2)d_r^3(1)}.
\end{equation}
 \end{theorem}
From Theorem \ref{thm2} we derive the main result of the paper.
\begin{theorem} \label{thm3}
Under assumptions of the Theorem \ref{thm2},
\begin{equation}\label{res3}
    \sqrt{k}\left( \hat{\gamma}^{(j)}_n(k,r)-\gamma \right) {\buildrel \rm d \over \to} \mathcal{N}\left( \mu \nu_j(r), \sigma^2_j(r)\right), \quad j=1,2,3,
\end{equation}
where
\begin{eqnarray*}
  \nu_1(r) &=& \frac{d_r(1)}{d_r(1) -\rho}, \quad \sigma_1^2 (r)=\frac{\gamma^2 d_r^2(1)}{d_r(2)},  \\
  \nu_2 (r) &=& \frac{d_r(1)(1-\rho-\gamma^2 r^2)}{(1+\gamma r) (d_r(1) -\rho)^2}, \quad \sigma_2^2 (r) = \frac{\gamma^2d_r^2(1)(d_r(2) +2\gamma^4 r^4)}{(1+\gamma r)^2 d_r^3(2)},  \\
  \nu_3 (r) &=& \frac{d_r^2(1)}{(d_r(1)-\rho)^2}, \quad \sigma_3^2 (r) =\frac{2\gamma^2d_r^4(1)}{d_r^3(2)}.
\end{eqnarray*}
\end{theorem}

Having asymptotic normality of the introduced estimators in Section 2 we compare the asymptotic mean square error (AMSE) of these estimators. In \cite{Paul16} we showed that the GHE $\hat{\gamma}^{(1)}_n(k,r)$ dominates the Hill estimator in all region of the  parameters $\g$ and $\rho$: $\g>0, \rho <0.$ Although this domination is rather small (see Fig. 2, the right graph ), but theoretically it is important, since, as far as we know, this is the first estimator, which asymptotically performs better than the Hill estimator in all region of possible range of two parameters $\g$ and $\rho$. Now we compare estimators $\hat{\gamma}^{(j)}_n(k,r), \ j=2,3,$  with  $\hat{\gamma}^{(1)}_n(k,r)$, and the results are promising. First of all, GMRE substantially outperforms MRE in all region of the  parameters  $\g>0, \rho <0,$ see Fig. 2, the left graph. Fig. 3, the right graph presents comparison of GMRE with GHE, and asymptotic result (solid line)  shows that no one estimator dominates in all region $-\infty <\rho <0$ (the ratio of AMSE does not depend on $\g$),
while simulation results (points) demonstrates the domination of GMRE  for all values of $\rho$, for which the simulation was performed.

Finally, we formulate some results concerning the  robustness of the introduced estimators. We follow the paper \cite{Beran}, where robustness was considered for the HME (or in our notation, for $\hat{\gamma}^{(1)}_n(k,r)$).

  To define the robustness measure for the estimator $\hat{\gamma}^{(\ell)}_n(k,r)$, instead of
$\hat{\gamma}^{(\ell)}_n(k,r)$ we will use the notation $\hat{\gamma}^{(\ell)}_n(k,r; X_1, \dots, X_n)$ .
Let us define
$$
\Delta \hat{\gamma}^{(\ell)}_n(k,r,x)=\hat{\gamma}^{(\ell)}_n(k,r; X_1, \dots, X_{n-1},x)- \hat{\gamma}^{(\ell)}_{n-1}(k-1,r; X_1, \dots, X_{n-1}).
$$
Then, following \cite{Beran}, for fixed $n$ and $k$  we take a quantity
$$
B_n^{(\ell)} (k,r)=\lim_{x \to \infty} \Delta \hat{\gamma}^{(\ell)}_n(k,r,x)
$$
which measures the worst effect of one arbitrarily large
contamination on the estimator $\hat{\gamma}^{(\ell)}_n(k,r)$. For fixed $n$ and $k$ these quantities are random variables, but it turns out that asymptotically they become constants, depending on $\g$ and $r$ (here it is appropriate to note, that results on robustness are based on Theorem \ref{thm1}, thus there is no dependence on $\rho$).

Let us denote by $B^{(\ell)}:=B^{(\ell)} (\g,r)$ the limit in probability  of $B_n^{(\ell)} (k,r)$, $\ell=1,2,3$, as $n\to \infty$ and (\ref{cond3}) holds.
\begin{theorem} \label{thm4}
Suppose that $X_1,\dots, X_n$ are i.i.d. nonnegative random
variables whose quantile function $U$ satisfies condition
(\ref{cond2}). Let $\g r<1$.
Then we have
$$
B^{(j)} =
\left\{
  \begin{array}{ll}
    0, & \hbox{$r<0$,} \\
    \infty, & \hbox{$r=0$,} \\
   (1-\gamma r)/r, & \hbox{$0<r<1/\gamma$.}
  \end{array}
\right. \quad j=1,2,3.
$$
\end{theorem}

Assuming the SORV condition (\ref{cond4}) we were able to find optimal values of $r$ for $\hat{\gamma}^{(j)}_n(k,r)$, $j=1,3 $ (see formulas (\ref{comp6}) and (\ref{comp7}) in Section 2), therefore, for the generalized Hill estimator we get $B^{(1)}(\g, r_1^*)=\g (1- \rho+{\sqrt{(2- \rho)^2 -2}})$. For the generalized moment ratio estimator the situation is even better, since optimal value $r_3^*<0$, therefore $B^{(3)}(\g, r_3^*)=0$,
while $B^{(j)}(\g, 0)=\infty$.
At first it seemed for us little bit strange, that for all three estimators we get the same value of $B^{(j)}$, but looking more carefully to the construction of this measure of robustness, we realized that it is quite natural and even the proof is almost trivial. Since the second term in the expression of   $\Delta \hat{\gamma}^{(\ell)}_n(k,r,x)$ is independent of $x$, moreover, it tends to $\g$, we have
\begin{equation} \label{robust1}
B_n^{(\ell)} (k,r)=\lim_{x \to \infty}\hat{\gamma}^{(\ell)}_n(k,r; X_1, \dots, X_{n-1},x)- \hat{\gamma}^{(\ell)}_{n-1}(k-1,r; X_1, \dots, X_{n-1}).
\end{equation}
Thus, robustness of the given estimator depends essentially on this first limit, which can be zero, infinity of finite,
depending on the term $g_{r,u}\left(x/X_{n-k,n} \right)$. For all classical estimators $\hat{\gamma}^{H}, \hat{\gamma}^{M}, \hat{\gamma}^{MR}$ ($r=0$)
this term tends to infinity, therefore we get infinite value for robustness, while for all three generalizations we are getting that this limit
as $x\to \infty$ is $1/r$, if $r>0$ and is $\g$, if $r<0$, thus we are getting result of Theorem \ref{thm4}. Moreover, the proof of this theorem shows that we can contaminate the sample
not by one large value, but by several, and the asymptotic result will be the same - generalized estimators will remain  (asymptotically) robust.

The rest of the paper is organized as follows. In the next section we investigate asymptotic mean square error of the introduced estimators, and compare these estimators with the generalized Hill estimator, using the same methodology as in \cite{Haan}, and provide some simulation results. Then there are formulated conclusions, and the last section is devoted to the proofs of the results. At the end of the proof of Theorem \ref{thm3} we discuss the alternative proof of this result based on the paper \cite{Drees}.

\section{Theoretical Comparison of the Estimators and Monte-Carlo simulations}

As in  \cite{Haan} we can write  the following relation  for the asymptotic mean squared error
of the estimator $\hat{\gamma}^{(\ell)}_n(k,r)$
\begin{equation} \label{comp01}
{\rm AMSE} \left(\hat{\gamma}^{(\ell)}_n(k,r)\right) \sim \nu_{\ell}^2(r) A^2\left( \frac{n}{k}\right)+\frac{\sigma_{\ell}^2(r)}{k}, \quad \ell=1,2,3,
\end{equation}
where a sequence $k=k(n)$ satisfies (\ref{cond5}). Here and below  we write  $a_n \sim b_n$ if $a_n / b_n \to 1$ as $n \to \infty$.
As  in \cite{Paul16}, assuming that $r$ is fixed, we perform the minimization of ${\rm AMSE} \left(\hat{\gamma}^{(\ell)}_n(k,r)\right)$
 with respect to $k$ . We will obtain optimal value of  $k(n)$ for the
estimator $\hat{\gamma}^{(\ell)}_n(k,r)$, which will be denoted by $k_{\ell}^* (r)$. Then we  minimize asymptotic mean squared error
${\rm AMSE} \left(\hat{\gamma}^{(\ell)}_n\left(k_{\ell}^* (r),r\right)\right)$ with respect to $r$.

Thus, let us assume that $r$ is fixed. We define the function $a$ by  the following relation:
\begin{equation}\label{a-def}
A^2(t) \sim \int_t^{\infty} a(u) \ \d u, \quad t \to \infty.
\end{equation}
 By applying Lemma 2.8 in \cite{DekkHaan}
we get that minimum of right hand side of (\ref{comp01}) is achieved with
\begin{equation}\label{comp02}
k_{\ell}^* (r) = \left( \frac{\sigma_{\ell}^2(r)}{\nu_{\ell}^2(r)}\right)^{1/(1-2\rho)} \frac{n}{a^{\leftarrow}(1/n)}.
\end{equation}
Following the lines in \cite{Haan}  and using  assumption (\ref{cond5}), we get
\begin{equation}\label{comp03}
k_{\ell}^* (r) A^2\left( \frac{n}{k_{\ell}^* (r)}\right) \sim \frac{\sigma_{\ell}^2(r)}{(-2\rho)\nu_{\ell}^2(r)}, \quad n \to \infty.
\end{equation}
Substituting (\ref{comp02}) and (\ref{comp03}) into (\ref{comp01}) we obtain
\begin{equation} \label{comp04}
{\rm AMSE} \left(\hat{\gamma}^{(\ell)}_n\left(k_{\ell}^* (r),r\right)\right) \sim \frac{1-2\rho}{(-2\rho)}
\left( \nu_{\ell}^2(r) \left(\sigma_{\ell}^2(r)\right)^{-2\rho}\right)^{1/(1-2\rho)} \frac{a^{\leftarrow}(1/n)}{n}.
\end{equation}
Now we minimize ${\rm AMSE} \left(\hat{\gamma}^{(\ell)}_n\left(k_{\ell}^* (r),r\right)\right)$ (more precisely, the right-hand side of (\ref{comp04})) with respect to $r$, and for this it is sufficient to
 minimize the product $\nu_{\ell}^2(r) \left(\sigma_{\ell}^2(r)\right)^{-2\rho}$ with respect to $r$.
Let us note that asymptotic parameters $\nu_{\ell} (r)$, $1\le \ell\le 3$ depend on parameters $\rho$ and $\gamma r$, while
quantities $\sigma_{\ell}^2(r)/\gamma^2$, $1\le \ell\le 3$ depend on the product $\gamma r$ only. Therefore, it is  convenient
to introduce notation $R = \gamma r$ and to consider minimization of the function
$$
\eta_{\ell}(R)=\gamma^{4 \rho} \nu_{\ell}^2(R/\gamma) \left(\sigma_{\ell}^2(R/\gamma)\right)^{-2\rho},
$$
with respect to $R$ satisfying inequality $R<1/2$. Equating the derivative of this function to zero, we get
\begin{equation}\label{comp5}
  \sigma_{\ell}^2(R/\gamma)\frac{\d  \nu_{\ell}(R/\gamma)}{\d R}-
\rho \nu_{\ell}(R/\gamma)   \frac{\d \sigma_{\ell}^2(R/\gamma)}{\d R}=0.
\end{equation}
By substituting the values of $\nu_{3}(R/\gamma)$ and  $\sigma_{3}^2(R/\gamma)$ into  equation (\ref{comp5})
we get the equation
$
R^2 -R(2-\rho)+\rho=0.
$
Whence it follows that the optimal value of $R$ for the
estimator $\hat{\gamma}^{(3)}_n\left(k_{3}^* (r),r\right)$ is
\begin{equation}\label{comp6}
R_3^*= \frac{(2-\rho)-\sqrt{(2-\rho)^2-4\rho}}{2},
\end{equation}
 since the second root of the quadratic equation does not satisfies the relation $R<1/2$. As for the estimator $\hat{\gamma}^{(1)}_n\left(k_{1}^* (r),r\right)$, the optimal value of the parameter $r$ was found in \cite{Paul16}, and in our notation (it is necessary to note, that SORV condition in \cite{Paul16} was used with different parametrization, see (1.3) therein) the optimal value of $R$ is
\begin{equation}\label{comp7}
R_1^*=\frac{(2-\rho)-\sqrt{(2-\rho)^2-2}}{2}.
\end{equation}
Unfortunately, the situation with the estimator $\hat{\gamma}^{(2)}_n\left(k_{2}^* (r),r\right)$ is more complicated.
By substituting  the expressions of  $\nu_{2}(R/\gamma)$ and  $\sigma_{2}^2(R/\gamma)$ into (\ref{comp5})
we get the  equation
\begin{eqnarray*}
  &&2R^9 -2R^8(1-\rho)-2R^7(5-3\rho)+2 R^6 (\rho^2-3\rho+6)-2\rho R^5(5-2\rho)-6 R^4(1-\rho)^2 \\
&& \ +R^3(8\rho^2-22\rho+15)-2R^2(5\rho^2-14\rho+9)+4R(\rho^2-3\rho+2)-(1-\rho)=0.
\end{eqnarray*}
For a fixed given value of $\rho$ this equation of $9$-th order was solved with "Wolfram Mathematica 6.0".
It turns out that depending on the parameter $\rho$ it has 3 real roots and three pairs of conjugate roots or
5 real roots and two pairs of conjugate roots. All real roots were substituted into the function $\eta_{2}(R)$ and optimal value was found in this way.
Numerical values of optimal value $R_2^*$  as a function of $\rho$ are provided in Fig.1. Although we got explicit and very simple expressions for  the optimal values $R_1^*$ and $R_3^*$, we provide these two functions (as functions of $\rho$) in the same  Fig.1.
\begin{figure}[ht!] \label{pic-a}
\begin{center}
 \includegraphics[width=0.60\textwidth]{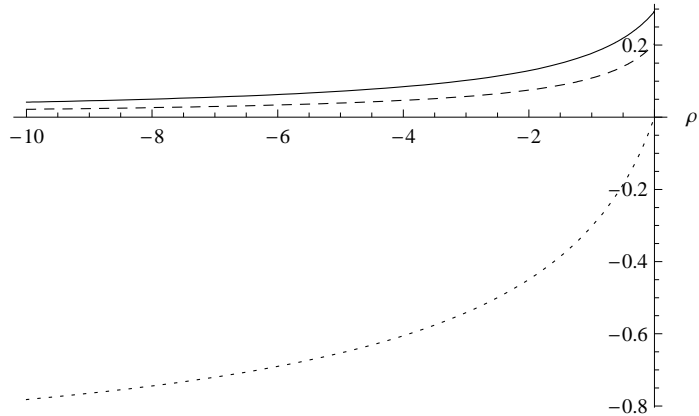}
\end{center}
\caption{Graph of the functions $R_1^* (\rho)$ (solid line), $R_2^*
(\rho)$ (dashed line) and $R_3^* (\rho)$ (dotted line),
$-10<\rho<0$}
\end{figure}

Let $r_{\ell}^* = R_{\ell}^* / \gamma$ denote optimal value of the parameter $r$.
From this picture we see that the first two  functions $R_i^*, \ i=1,2$ has comparatively small range of values (for $R_3^* (\rho)$ the range is $(-1, 0)$ ), this means that optimal value of parameters $r_1^*$ and $r_2^*$ are not sensitive to the parameter $\rho$, but more sensitive to $\g$.

Now we are able compare the generalized Hill estimator $\hat{\gamma}^{(1)}_n\left(k_{1}^* (r_1^*),r_1^*\right)$
with the another generalization of the Hill estimator $\hat{\gamma}^{(2)}_n\left(k_{2}^* (r_{2}^*),r_{2}^*\right)$ and
generalized moment ratio estimator $\hat{\gamma}^{(3)}_n\left(k_{3}^* (r_{3}^*),r_{3}^*\right)$.
But before performing comparison of these estimators, at first we demonstrate that
generalized moment ratio estimator $\hat{\gamma}^{(3)}_n\left(k_{3}^* (r_3^*),r_3^*\right)$
outperforms the initial moment ratio estimator $\hat{\gamma}^{(3)}_n\left(k_{3}^* (0),0\right)$ in the whole
area $\{(\gamma, \rho): \ \gamma>0, \ \rho<0\}$. Denoting
$$
\psi_{MR}(\rho)=\lim_{n \to \infty} \frac{{\rm AMSE} \left(\hat{\gamma}^{(3)}_n\left(k_{3}^* (0),0\right)\right)}
{{\rm AMSE} \left(\hat{\gamma}^{(3)}_n\left(k_{3}^* (r_3^*),r_3^*\right)\right)},
$$
it is not difficult to get that
$$
\psi_{MR}(\rho)= \left( \frac{2^{-8\rho} \left(v(\rho)-\rho \right)^4\left( v(\rho)-1+\rho\right)^{-6\rho}}{(1-\rho)^4 \left( v(\rho)+\rho\right)^{4-8\rho}}\right)^{1/(1-2\rho)},
$$
where $v(\rho)=\left((2-\rho)^2-4\rho \right)^{1/2}$. Since we must investigate this function on negative half-line $\{ \rho<0\}$, it is convenient to denote $-\rho=x$ and to write
$$
{\tilde \psi}_{MR} (x) =\psi_{MR} (-x) = (f(x))^{g(x)},
$$
with
$$
f(x)=  \frac{2^{-8x} \left({\tilde v}(x)+x \right)^4\left({\tilde v}(x)-1-x\right)^{6x}}{(1+x)^4 \left({\tilde v}(x)-x\right)^{4+8x}}, \ g(x)=\frac{1}{1+2x},
$$
where ${\tilde v}(x)=v(-x)$.
Taking logarithm of ${\tilde \psi}_{MR} (x)$,  using the fact that $f$ is product of several elementary functions, and using the simple  relation
$v(x)-x=4 +O(x^{-1})$, one can get
$$
\lim_{x\to \infty}\ln {\tilde \psi}_{MR} (x)=3\ln 3 -4\ln 2, \ {\rm or} \ \ \lim_{x\to \infty}{\tilde \psi}_{MR} (x)=\frac{27}{16}=1.6875
$$
 In a similar way one can get
$$
\lim_{x\to 0}\ln {\tilde \psi}_{MR} (x)=0 , \ {\rm or} \ \ \lim_{x\to 0} {\tilde \psi}_{MR} (x)=1.
$$
As a matter of fact, $\psi_{MR} (0)=1$, but considering asymptotic normality and AMSE of estimators under consideration we excluded the case $\rho=0,$
therefore we calculate this last limit. More difficult is to show that the function ${\tilde \psi}_{MR} (x)$   is
monotone and increasing (or $\psi_{MR} (\rho)$ is decreasing), we skip these considerations, only we mention that we use the fact that
the logarithmic derivative of a product of functions is a sum of logarithmic derivatives of these functions.
The graph of the function $ \psi_{MR}(\rho)$, $\rho<0$ is provided in Fig.2 in left.  In the same picture we gave also results (in form of separate points) of simulations, which will be explained below. Surprisingly, simulation results are even better than theoretical asymptotical result - most points are above the graph of $ \psi_{MR}(\rho)$.

 \begin{figure}[ht!] \label{pic-c}
\begin{center}
 \includegraphics[width=0.40\textwidth]{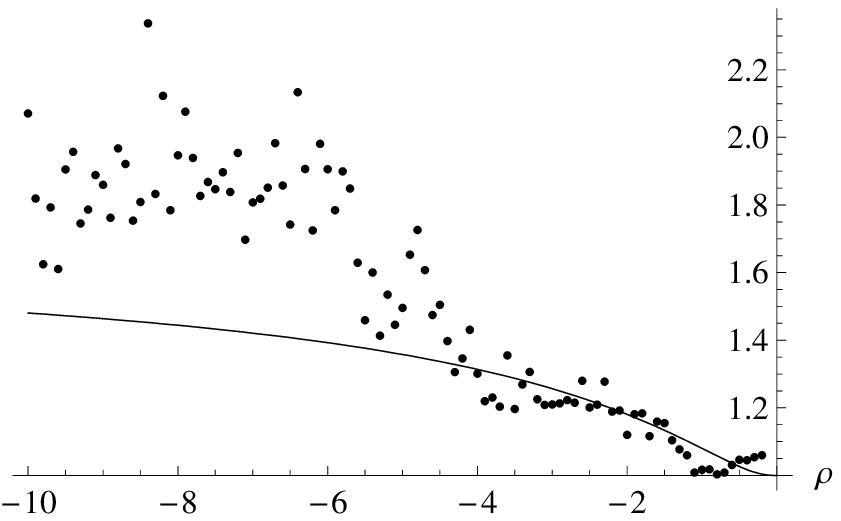} \ \includegraphics[width=0.40\textwidth]{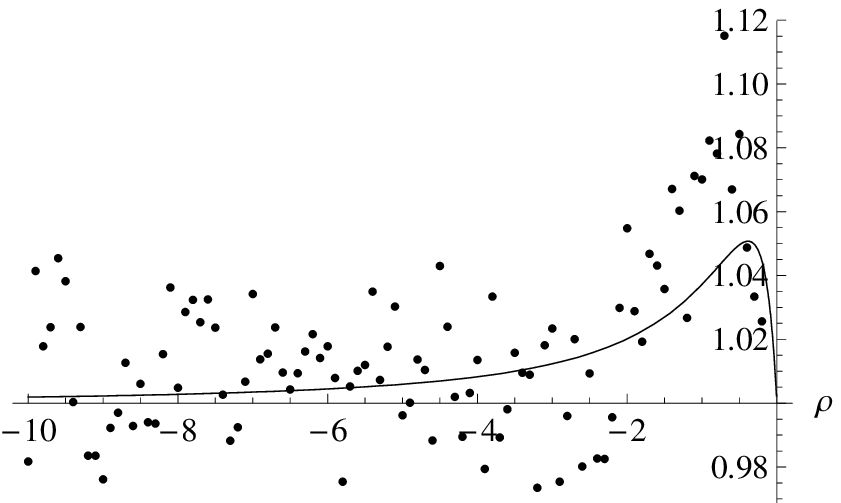}
\end{center}
\caption{Graph of the functions $\psi_{MR} (\rho)$
(solid line, in left),   $\psi_{H} (\rho)$ (solid line, in right) and results of Monte-Carlo simulations (points)}
\end{figure}

Now we can compare this result with the analogous  comparison of
the generalized  Hill estimator with the Hill estimator. If we
denote
$$
\psi_{H}(\rho)=\lim_{n \to \infty} \frac{{\rm AMSE}
\left(\hat{\gamma}^{(1)}_n\left(k_{1}^* (0),0\right)\right)} {{\rm
AMSE} \left(\hat{\gamma}^{(1)}_n\left(k_{1}^*
(r_1^*),r_1^*\right)\right)},
$$
then  we have
$$
\psi_{H}(\rho)=\left( \frac{(1-R_1^*
-\rho)^2(1-2R_1^*)^{-2\rho}}{(1-\rho)^2(1-R_1^*)^{2-4\rho}}\right)^{1/(1-2\rho)}.
$$
In \cite{Paul16}  it was shown that the improvement of the
generalized Hill estimator over the Hill estimator is important only
theoretically since the maximal value of the function $\psi_{H}$  is
$1.05$, see Fig.2 in right (this is the same graph as in Figure 2 in
\cite{Paul16}, only now we added results of simulation; in this case
simulation results are in good correspondence with theoretical
ones).
Thus, comparing with the generalized Hill estimator, we have
substantial improvement of the  moment ratio estimator, since in a
big range of the parameter $\rho$ the function $\psi_{MR}$ is bigger
than 1.3, and the maximal value is close to $1.7$.

Now we return to comparison of estimators $\hat{\gamma}^{(2)}_n\left(k_{2}^* (r_{2}^*),r_{2}^*\right)$ and
 $\hat{\gamma}^{(3)}_n\left(k_{3}^* (r_{3}^*),r_{3}^*\right)$ with the generalized Hill estimator $\hat{\gamma}^{(1)}_n\left(k_{1}^* (r_1^*),r_1^*\right)$, and we must investigate the following two functions:
$$
\varphi_{\ell} (\rho):=\lim_{n \to \infty} \frac{{\rm AMSE} \left(\hat{\gamma}^{(1)}_n\left(k_{1}^* (r_1^*),r_1^*\right)\right)}
{{\rm AMSE} \left(\hat{\gamma}^{(\ell)}_n\left(k_{\ell}^* (r_{\ell}^*),r_{\ell}^*\right)\right)}, \quad \ell=2,3.
$$
It is important to note that  both functions are independent of $\g$  and depend only on $\rho$.
In view of (\ref{comp04}) we have
$$
\varphi_{\ell} (\rho)=\left(
\frac{\nu_{1}^2\left(r_1^*\right) \left(\sigma_{1}^2\left(r_1^*\right)\right)^{-2\rho}}
{\nu_{\ell}^2\left(r_{\ell}^*\right) \left(\sigma_{\ell}^2\left(r_{\ell}^*\right)\right)^{-2\rho}}
\right)^{1/(1-2\rho)}, \quad \ell=2,3.
$$
Since we were able to obtain the optimal value of $R_2^*$ only numerically, we can provide only a numerically obtained graph of the function $\varphi_{2} (\rho)$, see Fig.3 on the left, therefore, in this case we did not provide simulation results.

\begin{figure}[ht!] \label{pic-b}
\begin{center}
 \includegraphics[width=0.40\textwidth]{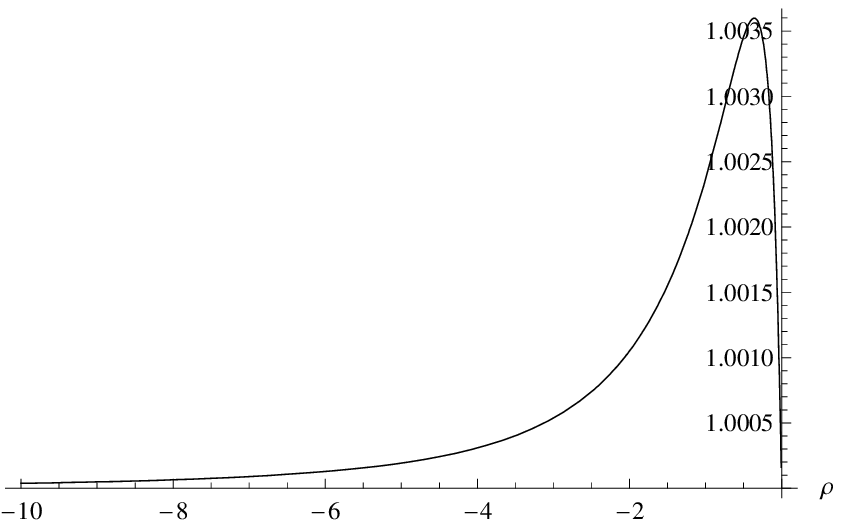} \  \includegraphics[width=0.40\textwidth]{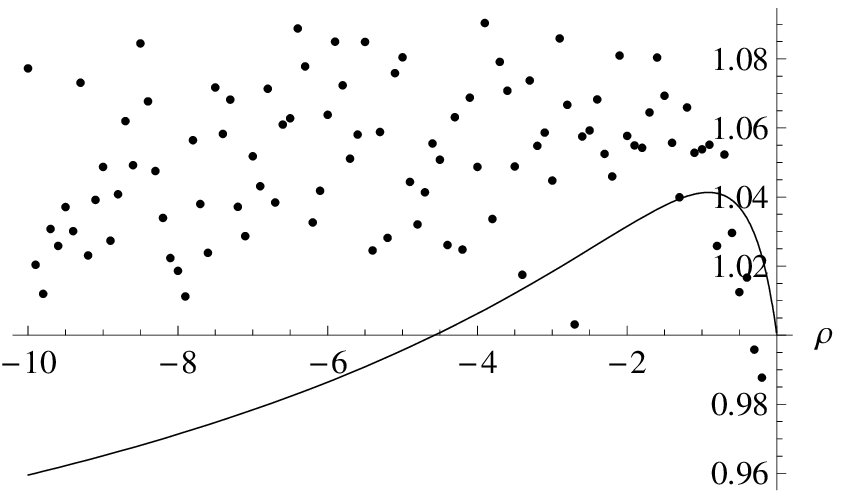}
\end{center}
\caption{Graph of the functions $\varphi_{2} (\rho)$ (on the left), $\varphi_{3} (\rho)$ (solid line, on the right) and
results of Monte-Carlo simulations (points)}
\end{figure}

Although this graph allow to believe that the new generalization of the Hill estimator dominates the generalized Hill estimator (which is the same as HME) in all region of parameters $\{(\gamma, \rho): \ \gamma>0, \ \rho<0\}$ , but without explicit expression of the function $\varphi_{2} (\rho)$ we can not prove this.

Finally, comparing two generalizations of the Hill and moment ratio estimators  (which, in our opinion are the most successful, since for both of them we have quite simple expression for optimal value of $R_{\ell}^*, \ell=1,3$ )we have
$$
\varphi_{3} (\rho)=\left(\frac{3^{-6 \rho} \left( v(\rho)-\rho\right)^{8-8\rho}\left( w(\rho)+(1-\rho)\right)^{-2\rho} }
{4^{3-5 \rho} (1-2\rho)^2 \left( v(\rho)+(1-\rho)\right)^{-6\rho} \left( w(\rho)-\rho\right)^{4-4\rho}} \right)^{1/(1-2\rho)},
$$
where $w(\rho)=\left((2-\rho)^2-2 \right)^{1/2}$. As can be seen
from Fig. 3 in right, generalized moment ratio estimator
$\hat{\gamma}^{(3)}_n\left(k_{3}^* (r_3^*),r_3^*\right)$ dominates
the generalized Hill estimator $\hat{\gamma}^{(1)}_n\left(k_{1}^*
(r_1^*),r_1^*\right)$  for $\rho \in (\tilde{\rho}, 0)$, where
$\tilde{\rho} \approx -4.57018$. It is not difficult to show that
$$
\lim_{\rho \to -\infty}\varphi_{3} (\rho)=\frac{3^3}{2^5}=0.84375.
$$
And empirical results, as in the case of function $\psi_{MR} (x)$ in
Fig. 2 (in left), are somehow unexpected - almost all empirical points in the
figure are not only  above the graph of $\psi_{MR} (x)$, but they are
bigger than 1, this means that the  empirical  MSE of generalized moment ratio
estimator is smaller than the empirical MSE of the generalized Hill estimator for
all values of $\rho$ in the interval $-10<\rho<0$, not only for
$\rho \in (\tilde{\rho}, 0)$, as gives asymptotic theoretical
result.

Now we shall provide some results of Monte-Carlo simulations (part of these results are given in Figures 2 and 3 together with theoretical results). We must admit that these results are very preliminary, since for practical application of the proposed estimators still we are facing with problems of estimating the parameters $\rho$ and $r$, and the main problem is with estimation of $\rho$ (and this problem is for almost all estimators, if we want to have asymptotic normality of the estimator), since  optimal value of $r$ is expressed as a function of $\rho$. We intend to return to this problem in a separate paper.

For simulations we shall use a little bit more restrictive condition than (\ref{cond4}), namely, we assume that  the distribution function
$F(x)$ under consideration belongs to the Hall's class of Pareto type distributions
(\cite{Hall},\cite{Hall1}), i.e.
\begin{equation}\label{n01}
1-F(x)=\left(\frac{x}{C}\right)^{-1/\gamma} \left(1+
\frac{\beta}{\rho} \left(\frac{x}{C}\right)^{\rho/\gamma} +o\left(
x^{\rho/\gamma}\right)\right), \quad x \to \infty,
\end{equation}
where $C>0$, $\beta \in \mathbf{R}\setminus \{0\}$ and $\rho<0$. This assumption is assumed for the following reason. Taking the ratio of AMSE of two estimators we do not need to know the function $a^{\leftarrow}$, but for simulations, having given sample size $n$, we must calculate the value of $k_{\ell}^* (r)$ in (\ref{comp02}) and the empirical MSE of estimators, and for this we must have the function $a^{\leftarrow}$. Assuming (\ref{n01}),
we have that the second order condition  (\ref{cond4}) holds with
$A(t)=\gamma \beta t^{\rho}$ and from (\ref{a-def}) it follows that
$a^{\leftarrow}(t)=\left( -2\rho \gamma^2
\beta^2\right)^{1/(1-2\rho)} t^{1/(2\rho-1)}$. Now we
can rewrite (\ref{comp02}) as follows:
$$
k_{\ell}^* (r,\beta, \rho) = \left( \frac{\sigma_{\ell}^2(r)}{-2\rho \beta^2 \gamma^2 \nu_{\ell}^2(r)}\right)^{1/(1-2\rho)} n^{-2\rho/(1-2\rho)}.
$$
In fact, quantities $k_{\ell}^* (0,\beta, \rho)$, $\ell=1,3$  depend on $\beta$, $\rho$ and $n$ only, thus replacing
$\beta$ and $\rho$ by  some  estimators $\hat{\beta}_n$ and $\hat{\rho}_n$, we obtain the empirical values of the parameter $k(n)$
for the Hill and the moment ratio estimators:
\begin{eqnarray} \label{k1-H}
\hat{k}_{1,n} &=& \left( \frac{(1-\hat{\rho}_n)^2}{-2\hat{\rho}_n  \hat{\beta}_n^2}\right)^{1/(1-2\hat{\rho}_n)} n^{-2\hat{\rho}_n/(1-2\hat{\rho}_n)}, \\
\hat{k}_{3,n} &=& \left( \frac{(1-\hat{\rho}_n)^4}{-\hat{\rho}_n \hat{\beta}_n^2}\right)^{1/(1-2\hat{\rho}_n)} n^{-2\hat{\rho}_n/(1-2\hat{\rho}_n)}.
\end{eqnarray}
For corresponding generalized estimators we have additionally to take estimators of optimal parameter $R,$ therefore we have
\begin{eqnarray}
\label{k1-GH}
\hat{K}_{1,n} &=& \left( \frac{(1-\hat{\rho}_n-R_1^*\left(\hat{\rho}_n\right))^2}{-2\hat{\rho}_n  \hat{\beta}_n^2 (1-2R_1^*\left(\hat{\rho}_n\right))}\right)^{1/(1-2\hat{\rho}_n)}
n^{-2\hat{\rho}_n/(1-2\hat{\rho}_n)}, \\
\hat{K}_{3,n} &=& \left( \frac{(1-\hat{\rho}_n-R_3^*\left(\hat{\rho}_n\right))^4}{-\hat{\rho}_n   \hat{\beta}_n^2(1-2R_3^*\left(\hat{\rho}_n\right))^3}\right)^{1/(1-2\hat{\rho}_n)} n^{-2\hat{\rho}_n/(1-2\hat{\rho}_n)}.
\end{eqnarray}
Thus, in our simulations the comparison is made between the Hill estimator
$\hat{\gamma}^{(1)}_n\left(\hat{k}_{1,n},0\right)$, the generalized Hill estimator
$\hat{\gamma}^{(1)}_n\left(\hat{K}_{1,n},r_1^*\left(\hat{\gamma}^{(1)}_n\left(\hat{k}_{1,n},0\right), \hat{\rho}_n\right)\right)$, the moment ratio estimator
$\hat{\gamma}^{(3)}_n\left(\hat{k}_{3,n},0\right)$ and the generalized moment ratio estimator
$\hat{\gamma}^{(3)}_n\left(\hat{K}_{3,n},r_3^*\left( \hat{\gamma}^{(3)}_n\left(\hat{k}_{3,n},0\right), \hat{\rho}_n\right)\right)$, with parameters given in (\ref{k1-H})-(2.13).

We generated 1000 times samples  $X_1,\dots ,X_n$ of i.i.d. random variables of  size $n = 1000$ with  the following two d. f. with the extreme value
index $\gamma$, parameter $\rho$ and satisfying (\ref{n01}):

(i) the Burr  d.f. $F(x) = 1 - \left(1 + x^{-\rho/\gamma}\right)^{1/\rho}$, $x \ge 0$;

(ii) the Kumaraswamy generalized exponential  d.f. $F(x)=1-\left( 1-\exp\left\{ -x^{\rho/\gamma}\right\}\right)^{-1/\rho}$, $x \ge 0$.

 The parameter $\beta$ which is present in (\ref{n01}) for the Burr distribution is $1$ and for the Kumaraswamy distribution - $1/2$, and $C=1$ for both distributions.
To calculate the Hill and the generalized Hill estimators we used the following algorithm:

1.  Estimate the parameter $\rho$ by the following estimator proposed in \cite{Fraga1}:
$$
\hat{\rho}_n (k,\tau)=-\left| 3\left(T_n^{(\tau)} (k)-1 \right)\left (T_n^{(\tau)} (k)-3\right )\right|,
$$
where
$$
T_n^{(\tau)} (k)= \left (\left(G_n(k, 0,1)\right)^{\tau}-\left(G_n(k, 0,2)/2\right)^{\tau/2}\right )
\left ({\left(G_n(k, 0,2)/2\right)^{\tau/2}-\left(G_n(k, 0,3)/6\right)^{\tau/3}}\right )^{-1}
$$
with $\tau>0$, and
$$
T_n^{(0)} (k)=\frac{\ln\left(G_n(k, 0,1) \right)-(1/2)\ln\left(G_n(k, 0,2)/2 \right) }
{(1/2)\ln\left(G_n(k, 0,2)/2 \right)-(1/3)\ln\left(G_n(k, 0,3)/6 \right)}.
$$
To decide which values of parameters $\tau$ ($0$ or $1$) and $k$ to take in the above written estimator  we realized the algorithm provided in \cite{Car3}.

2. To estimate  the parameter $\beta$ use the estimator $\hat{\beta}_n \left(k, \hat{\rho}_n (k,\tau)\right)$, where
$$
\hat{\beta}_n \left(k, \rho\right)=\frac{\left(\frac{k}{n}\right)^{\rho} \left\{\left(\frac{1}{k} \sum_{i=1}^k \left(\frac{i}{k}\right)^{-\rho} \right)\left(\frac{1}{k} \sum_{i=1}^k W_i \right)-
\left(\frac{1}{k} \sum_{i=1}^k \left(\frac{i}{k}\right)^{-\rho} W_i\right)\right\}}{\left(\frac{1}{k} \sum_{i=1}^k \left(\frac{i}{k}\right)^{-\rho} \right)\left(\frac{1}{k} \sum_{i=1}^k \left(\frac{i}{k}\right)^{-\rho} W_i\right)
-\left(\frac{1}{k} \sum_{i=1}^k \left(\frac{i}{k}\right)^{-2\rho} W_i\right)}
$$
and  $W_i=i \ln \left (X_{n+i-1,n}/ X_{n+i,n}\right )$, $1 \le i \le k <n$. This estimator was introduced in \cite{Gomes-a}. Also, as in the step 1 to choose the parameter $k$ we applied algorithm from \cite{Car3}.


3. By using (\ref{k1-H}) estimate parameter $k$ for the Hill estimator and then obtain $\hat{\gamma}^{(1)}_n\left(\hat{k}_{1,n},0\right)$;

4. Estimate $R_1^*\left(\hat{\rho}_n\right)$ and $r_1^*\left(\hat{\gamma}^{(1)}_n\left(\hat{k}_{1,n},0\right), \hat{\rho}_n\right)$;

5. By using (\ref{k1-GH}) estimate parameter $k$ for the generalized Hill estimator and find
estimate $\hat{\gamma}^{(1)}_n\left(\hat{K}_{1,n},r_1^*\left(\hat{\gamma}^{(1)}_n\left(\hat{k}_{1,n},0\right), \hat{\rho}_n\right)\right)$.

\noindent We used the similar algorithm (with obvious changes) for the moment ratio estimator $\hat{\gamma}^{(3)}_n\left(\hat{k}_{3,n},0\right)$ and the generalized moment ratio estimator $\hat{\gamma}^{(3)}_n\left(\hat{K}_{3,n},r_3^*\left( \hat{\gamma}^{(3)}_n\left(\hat{k}_{3,n},0\right), \hat{\rho}_n\right)\right)$. Having values of the estimators we calculate MSE and bias of these estimators and also the ratios of some pairs of MSE.

The results of simulations are summarized in Fig. 4 (for the Burr distribution) and in Fig. 5 (for the Kumaraswamy distribution).
For  the Burr distribution we took parameters $\gamma$ and $\rho$ in the intervals $(0,4]$ and $(-5,-0.2],$ respectively. In Fig. 4 (left)  we divided
this rectangle $(0,4] \times (-5,-0.2]$ into
squares $\Delta_{i,j}=(i/10, (i+1)/10] \times (-(j+1)/10, -j/10]$, $i=0,1,\dots,40$, $j=2,3,\dots,50$. By taking true values of parameters $\gamma$ and $\rho$
as coordinates of the center
 of the rectangle $\Delta_{i,j}$, we performed Monte-Carlo simulations. We colored the square $\Delta_{i,j}$ in black if
empirical MSE of the estimator $\hat{\gamma}^{(1)}_n\left(\hat{K}_{1,n},r_1^*\left(\hat{\gamma}^{(1)}_n\left(\hat{k}_{1,n},0\right), \hat{\rho}_n\right)\right)$
is the smallest among all four estimators under consideration, while
areas of domination of the estimators $\hat{\gamma}^{(1)}_n\left(\hat{k}_{1,n},0\right)$,
$\hat{\gamma}^{(3)}_n\left(\hat{k}_{3,n},0\right)$ and
$\hat{\gamma}^{(3)}_n\left(\hat{K}_{3,n},r_3^*\left( \hat{\gamma}^{(3)}_n\left(\hat{k}_{3,n},0\right), \hat{\rho}_n\right)\right)$
are in dark grey, grey and in white, respectively. In Fig. 4 (right) there are given
the  areas of domination of absolute value of the bias using the same colors as in the left picture.

\begin{figure}[ht!] \label{pic-f1}
\begin{center}
 \includegraphics[width=0.40\textwidth]{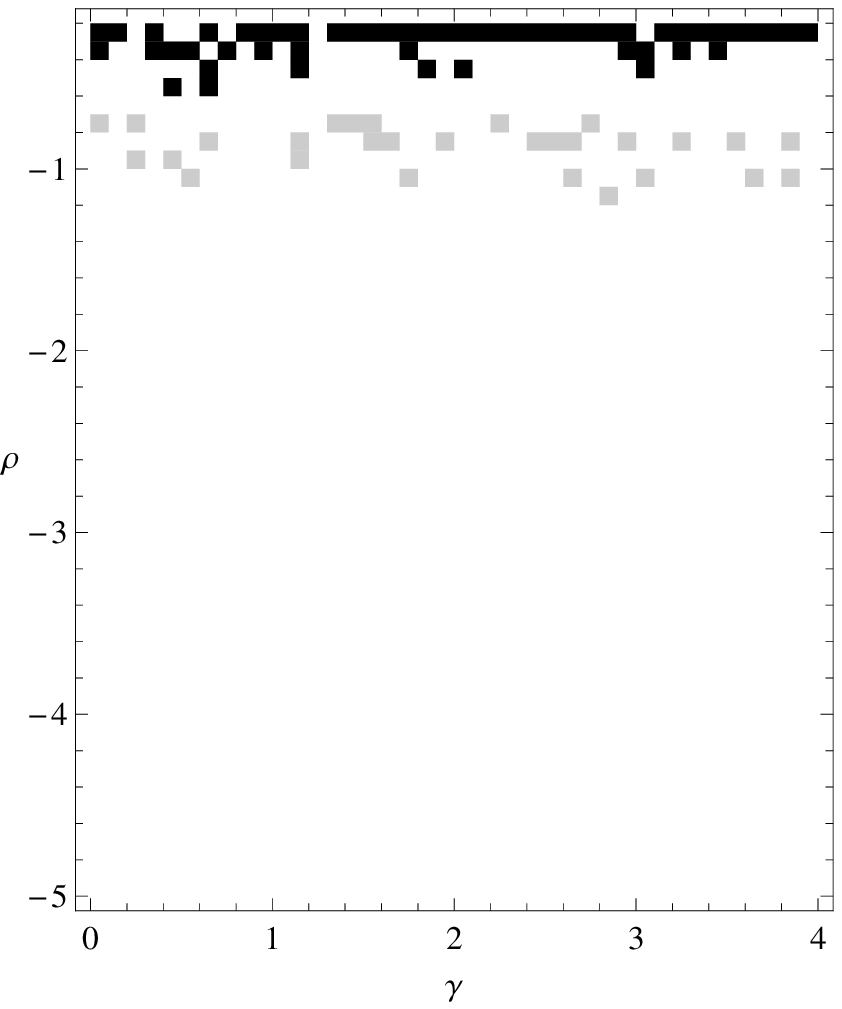} \  \includegraphics[width=0.40\textwidth]{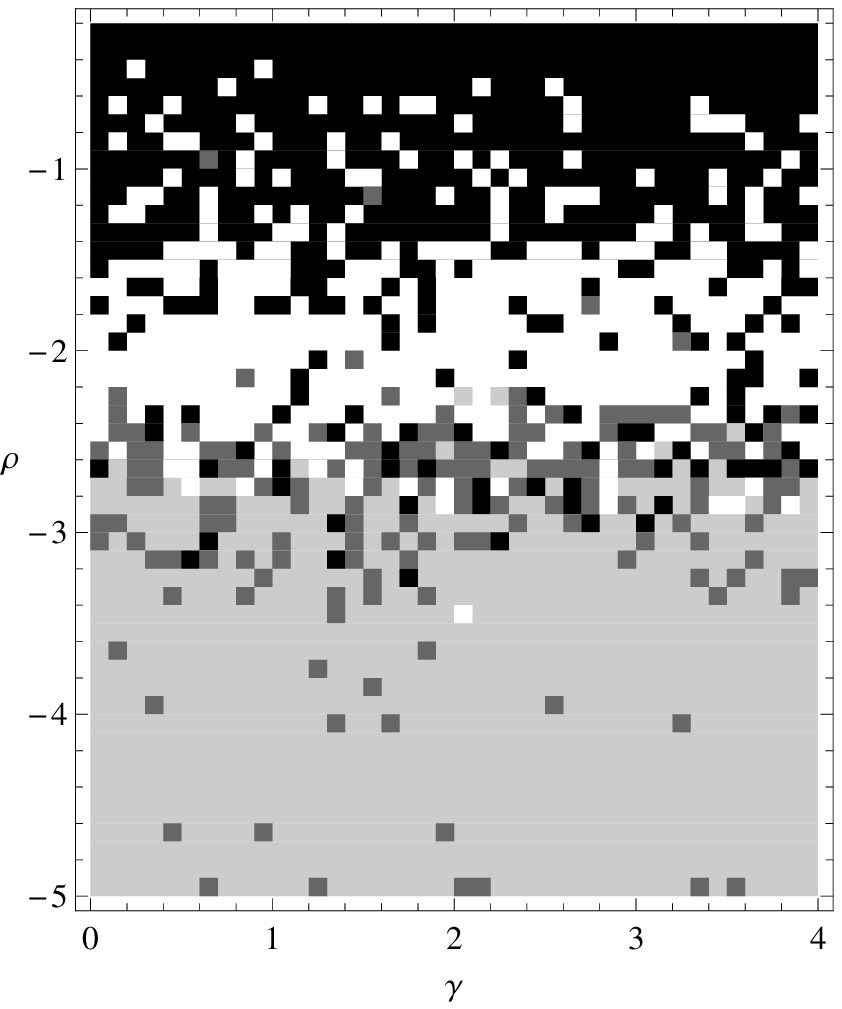}
\end{center}
\caption{Empirical comparison of the estimators by using Burr distribution }
\end{figure}

We faced with serious difficulties  while estimating  the second order parameter $\rho$
when  samples $X_1,\dots ,X_n$ are  simulated from Kumaraswamy distribution with $\rho<-2$ (for many samples computer was getting infinite value of sums
involved in the estimator of $\rho$), thus we restrict our simulations
in the rectangle $(\gamma, \rho) \in (0,4] \times [-2,-0.2]$. The corresponding areas of domination of MSE and bias for this distribution are presented in Fig.5.
\begin{figure}[ht!] \label{pic-f}
\begin{center}
 \includegraphics[width=0.40\textwidth]{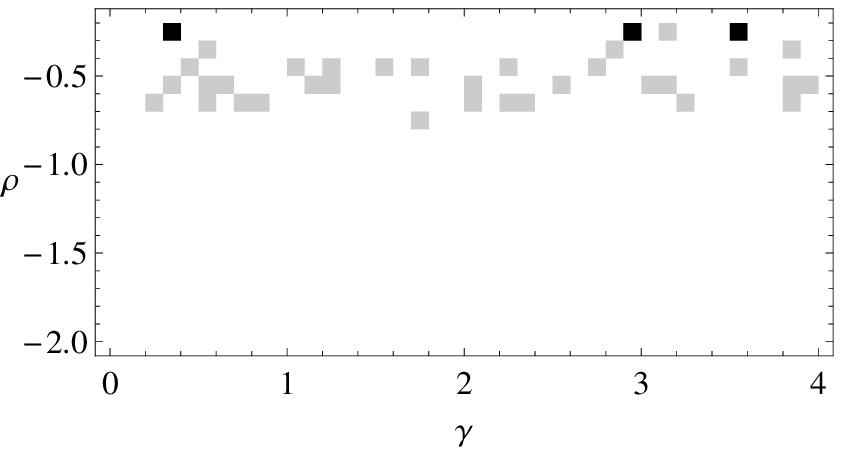} \  \includegraphics[width=0.40\textwidth]{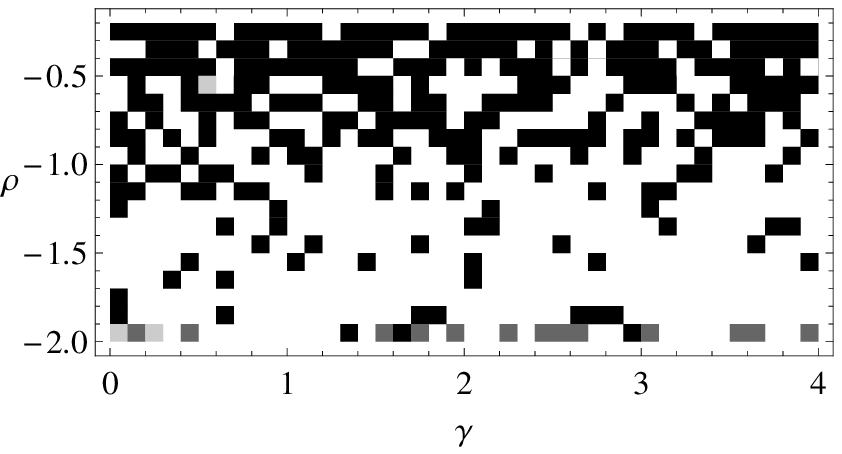}
\end{center}
\caption{Empirical comparison of the estimators by using Kumaraswamy distribution}
\end{figure}

Fig. 4 and 5  demonstrates that areas of domination almost do not depend on parameter $\gamma$ and essentially depend only on $\rho$.
This corresponds well to theoretical results that functions $\varphi_{\ell}(\rho)$, $\ell=2,3$,
$\psi_H(\rho)$, and $\psi_{MR}(\rho)$ depend on $\rho$ only.
Therefore, taking the particular value  $\gamma=1$ and the hundred of values of $\rho$ in the interval $(-10, 0)$ for the
Burr distribution we obtained ratios of empirical MSE's to complement theoretical comparison
and included these ratios (as separate  points)  in Fig. 2 and 3.
From Fig. 2 (right) we see that empirical results on comparison of MSE of the generalized Hill estimator with MSE of the Hill
estimator reflects quite well theoretical (asymptotic) function $\psi_{H} (\rho)$ - points are on both sides of this function.
Little bit unexpectedly, Fig. 2 (left) and 3 (right) reveals that empirical results differ from theoretical ones.
In Fig. 2 (left), where generalized moment ratio estimator is compared with the moment ratio estimator, in the interval $-5 \le \rho <0$
empirical points are very close to the theoretical function  $\psi_{H} (\rho)$, but  in the interval $-10<\rho <-5$ all points are above the theoretical curve,
this  means that the moment ratio estimator has MSE almost two times bigger that MSE of the generalized moment ratio estimator.
Empirical results in Fig. 3 (right) shows that the generalized moment ratio estimator performs better than GHE for all values of $\rho$ in the
interval $-10<\rho < 0$, while theoretical result predict such result only for $\rho \in (\tilde{\rho}, 0)$.

\section{Conclusions}
 We introduced a new parametric class of functions $g_{r,u}$ which allows to construct many new generalizations of well-known estimators, including such as the Hill, the moment, and the moment ratio estimators. We prove asymptotic normality of all these generalized estimators in a unified way and demonstrate that in the sense of AMSE new estimators are better than the classical ones, especially promising looks GMR estimator. We hope that this new parametric class of functions will be useful in the difficult problem of estimating the second order parameter $\rho$.

 Preliminary simulation results show quite good correspondence with the obtained theoretical results, but we admit that future work on the behavior of new estimators for middle size samples is needed.

\section{Proofs}

\noindent {\it Proof of the Theorem \ref{thm1}.} There is nothing to prove in the trivial case $r=u=0$. Keeping in mind
relation (\ref{sarysis}),  conclusion (\ref{res1}) is the immediate consequence of the Theorem 1 in \cite{Paul16} in the case $u=0$ and $\gamma r <1$. The case $r=0, \ u \ge 1$ was investigated in \cite{GomesMar}.

Consider the case $r \not =0$ and $u>-1$, $u \not =0$. Let us recall that the
function $U(t), \ t \ge 1$ varies regularly at infinity with the
index $\gamma$, thus, by applying Potter's bound we have: for
arbitrary $\epsilon>0$ there exits $t_0$, such that, for $ x \ge 1$
and $t \ge t_0$,
\begin{equation}\label{Pot_in}
\left((1-\epsilon)x\right)^{\gamma-\epsilon} < \frac{U(tx)}{U(t)} <
\left((1+\epsilon)x\right)^{\gamma+\epsilon}.
\end{equation}
In order to apply (\ref{Pot_in}) for the function $t^r$ it is
convenient to introduce the notation $\epsilon_{\pm}(r)$, where
$\epsilon_{\pm}(r)=\epsilon$, if $r>0$, and $\epsilon_{\pm}(r)=-\epsilon$,
if $r<0$. Then we get
\begin{eqnarray} \label{Pot_r}
\left\{(1-\epsilon_{\pm}(r))x\right\}^{r(\gamma-\epsilon_{\pm}(r))}  <
\left(\frac{U(tx)}{U(t)}\right)^r <
\left\{(1+\epsilon_{\pm}(r))x\right\}^{r(\gamma+\epsilon_{\pm}(r))}.
\end{eqnarray}
Similarly we get  the following inequalities
\begin{equation}\label{Pot_log}
(\gamma-\epsilon_{\pm}(u))^u \ln^u\left((1-\epsilon_{\pm}(u))x\right) < \ln^u \left(
\frac{U(tx)}{U(t)}\right) < (\gamma+\epsilon_{\pm}(u))^u \ln^u
\left((1+\epsilon_{\pm}(u))x\right).
\end{equation}
By multiplying inequalities (\ref{Pot_r}) and (\ref{Pot_log}), we
obtain
\begin{eqnarray*}
 c_1 x^{r(\gamma-\epsilon_{\pm}(r))}  \ln^u\left((1-\epsilon_{\pm}(u))x\right)
 < g_{r,u} \left( \frac{U(tx)}{U(t)}\right) <
  c_2 x^{r(\gamma+\epsilon_{\pm}(r))} \ln^u \left((1+\epsilon_{\pm}(u))x\right),
\end{eqnarray*}
where
\begin{eqnarray*}
  c_1 &=&  (\gamma-\epsilon_{\pm}(u))^u (1-\epsilon_{\pm}(r))^{r(\gamma-\epsilon_{\pm}(r))},  \\
  c_2 &=& (\gamma+\epsilon_{\pm}(u))^u (1+\epsilon_{\pm}(r))^{r(\gamma+\epsilon_{\pm}(r))}.
\end{eqnarray*}
Let $Y_1, \dots, Y_n$ be i.i.d. random variables with distribution
function $G(x)=1-(1/x)$, $x \ge 1$. Taking
\begin{equation}\label{pr01-0}
    t=Y_{n-k,n}, \quad x=Y_{n-i,n}/Y_{n-k,n}
\end{equation}
for $i=0,1,\dots, k-1$ we get
\begin{eqnarray} \label{pr01-0a}
c_1
\left(\frac{Y_{n-i,n}}{Y_{n-k,n}}\right)^{r(\gamma-\epsilon_{\pm}(r))}
\ln^u\left((1-\epsilon_{\pm}(u)) \frac{Y_{n-i,n}}{Y_{n-k,n}}\right)
 < g_{r,u} \left( \frac{U(Y_{n-i,n})}{U(Y_{n-k,n})}\right)  \\
\label{pr01-0b}
 c_2 \left(\frac{Y_{n-i,n}}{Y_{n-k,n}}\right)^{r(\gamma+\epsilon_{\pm}(r))} \ln^u \left((1+\epsilon_{\pm}(u))\frac{Y_{n-i,n}}{Y_{n-k,n}}\right)
 > g_{r,u} \left( \frac{U(Y_{n-i,n})}{U(Y_{n-k,n})}\right).
\end{eqnarray}
Note that $U(Y_i)=X_i$, $i=1,2, \dots,n$, thus
\begin{equation}\label{pr01-1}
G_n(k,r,u)=\frac{1}{k} \sum_{i=0}^{k-1} g_{r,u} \left(
\frac{U(Y_{n-i,n})}{U(Y_{n-k,n})}\right).
\end{equation}
From this equality, by summing inequalities (\ref{pr01-0a})
and(\ref{pr01-0b}), we get
\begin{eqnarray} \label{pr02-a}
\frac{c_1}{k} \sum_{i=0}^{k-1}
\left(\frac{Y_{n-i,n}}{Y_{n-k,n}}\right)^{r(\gamma-\epsilon_{\pm}(r))}
\ln^u\left((1-\epsilon_{\pm}(u)) \frac{Y_{n-i,n}}{Y_{n-k,n}}\right) <
G_n(k,r,u),\\
\label{pr02-b} \frac{c_2}{k} \sum_{i=0}^{k-1}
\left(\frac{Y_{n-i,n}}{Y_{n-k,n}}\right)^{r(\gamma+\epsilon_{\pm}(r))}
\ln^u \left((1+\epsilon_{\pm}(u))\frac{Y_{n-i,n}}{Y_{n-k,n}}\right)
>G_n(k,r,u).
\end{eqnarray}
By means of the  standard argument (see e.q. \cite{Haanbook}) one can
deduce that the left-hand side of (\ref{pr02-a}) equals (in
distribution) to the sum
\begin{eqnarray*}
\frac{c_1}{k} \sum_{i=1}^{k} Y_i^{r(\gamma-\epsilon_{\pm}(r))} \left(\ln(Y_i)+\ln(1-\epsilon_{\pm}(u))\right)^u.
\end{eqnarray*}
The expectation of this quantity  equals to $c_1\left(\Delta_1+ \Delta_2\right)$, where
\begin{eqnarray*}
\Delta_1&=& \int_1^{\infty} x^{r(\gamma -\epsilon_{\pm}(r))-2} \ln^u(x) \d x, \\
\Delta_2 &=& \int_1^{\infty} x^{r(\gamma -\epsilon_{\pm}(r))-2} \ln^u(x)\left\{\left(1+\frac{\ln(1-\epsilon_{\pm}(u))}{\ln(x)}\right)^u-1 \right\}\d x.
\end{eqnarray*}

One can verify that
\begin{equation}\label{id01}
    \int_1^{\infty} x^{a-2} \ln^b(x) \d x=(1-a)^{-1-b} \Gamma (1+b), \quad a<1, \ b>-1.
\end{equation}

By using the last identity,  assumptions
$1-\gamma r>0$, $u>-1$ and the fact that $r\epsilon_{\pm}>0$ we get
$$
c_1 \Delta_1 =
\frac{(1-\epsilon_{\pm}(r))^{r(\gamma-\epsilon_{\pm}(r))}
(\gamma-\epsilon_{\pm}(u))^u \Gamma(u+1)}{(1-r(\gamma-\epsilon_{\pm}(r)))^{u+1}} .
$$
Whence we get  $c_1 \Delta_1 \to \xi(r,u)$, as $\epsilon \to 0$.

Consider a quantity $\Delta_2$ now. If $0<u \le 1$, then $\epsilon_{\pm}(u)=\epsilon$ and
we use the following inequality which holds for any real numbers $a$ and $b$: $\left| |a|^u- |b|^u \right| \le |a-b|^u$. We have
\begin{eqnarray*}
\left|\Delta_2\right| &\le & \left|\ln(1-\epsilon)\right|^u \int_1^{\infty} x^{r(\gamma -\epsilon_{\pm}(r))-2}  \d x \\
&=& \frac{\left|\ln(1-\epsilon)\right|^u}{1-r(\gamma-\epsilon_{\pm}(r))}.
\end{eqnarray*}
Now it follows that $c_1 \Delta_2 \to 0$, as $\epsilon \to 0$. In the case $-1<u<0$ we have $\epsilon_{\pm}(u)=-\epsilon$, and since $0<-u<1$, applying the same inequality as above, we get
\begin{eqnarray*}
\left|\left(\ln(x)+\ln(1+\epsilon)\right)^u-  \ln^u(x)\right|&=&  \frac{\left|\ln(x)^{-u}-\left(\ln(x)+\ln(1+\epsilon)\right)^{-u}\right|}{\ln(x)^{-u}\left(\ln(x)+\ln(1+\epsilon)\right)^{-u}}\\
&\le& \frac{\left|\ln(1+\epsilon)\right|^{-u}}{\ln(x)^{-u}\left(\ln(x)+\ln(1+\epsilon)\right)^{-u}}.
\end{eqnarray*}
Now we take small $\delta>0$ such that $u-\delta>-1$, for example, one can take $\delta=1/2(u+1)$.
Keeping in mind that  $\ln(1+\epsilon)>0$, we can estimate
$$
\frac{1}{\left(\ln(x)+\ln(1+\epsilon)\right)^{-u}}\le \frac{1}{\left(\ln(x)\right )^\delta \left (\ln(1+\epsilon)\right)^{-u-\delta}}.
$$
Collecting the last two estimates we get
\begin{eqnarray*}
\left|\Delta_2\right| &\le & \left|\ln(1+\epsilon)\right|^{\delta} \int_1^{\infty} x^{r(\gamma -\epsilon_{\pm}(r))-2} \ln^{u-\delta}(x) \d x \\
&=& \frac{\left(\ln(1+\epsilon)\right)^\delta \Gamma(1+u-\delta)}{\left(1-r(\gamma-\epsilon_{\pm}(r))\right)^{1+u-\delta}}.
\end{eqnarray*}
This allows to deduce that $c_1 \Delta_2 \to 0$, as $\epsilon \to 0$.

Let $u>1$. By using inequality $\left| |a|^u- |b|^u \right| \le u \left( |a|^{u-1}- |b|^{u-1}\right)|a-b|$, which holds for
any real numbers $a$ and $b$, we get
\begin{eqnarray*}
\left|\Delta_2\right| &\le& 2 u \left| \ln(1-\epsilon_{\pm}(u)) \right|\int_1^{\infty} x^{r(\gamma -\epsilon_{\pm}(r))-2} \ln^{u-1}(x) \d x \\
&=& \frac{2 u \Gamma(u) \left| \ln(1-\epsilon_{\pm}(u)) \right|}{\left( 1-r(\gamma-\epsilon_{\pm}(r)) \right)^u},
\end{eqnarray*}
which implies $c_1 \Delta_2 \to 0$, $\epsilon \to 0$.

Thus, by applying the Khintchine weak law of large numbers, left-hand side of (\ref{pr02-a}) converges to zero in probability.
In a  similar way we can prove that the left-hand side of (\ref{pr02-b}) tends to zero in probability, too. Theorem \ref{thm1} is proved. \halmos

\bigskip

\noindent {\it Proof of Corollary \ref{cor1}.} From (\ref{res1}) it follows
$$
\left( G_n(k,r_1,u_1), \dots, G_n(k,r_{\ell},u_{\ell}) \right) {\buildrel \rm P \over \to}
\left( \xi(r_1, u_1), \dots, \xi(r_{\ell}, u_{\ell})
\right),
$$
as $n \to \infty$. Applying    Theorem 2.7 in \cite{Billingsley} we derive (\ref{cor1-res}).
\halmos

\bigskip

\noindent {\it Proof of the Theorem \ref{thm2}.} Let $r=0$. Theorem 3.2.5 in \cite{Haanbook} states
\begin{equation}\label{pthm2-01}
    \sqrt{k} \left(G_n(k,0,1) -\gamma\right) {\buildrel \rm d \over \to} \mathcal{N}\left( \frac{\mu}{1-\rho}, \gamma^2 \right), \quad n \to \infty.
\end{equation}
 The relation (\ref{pthm2-01}), together with $\sqrt{k} \left(G_n(k,0,0) -1\right) {\buildrel \rm P \over \to}  \ 0$ and  Theorem 3.9 in \cite{Billingsley}, gives (\ref{res2}) for $r=0$.

Consider now the case $r\not = 0$. Adjusting Potter's type bounds (3.4) in  \cite{Paul16} for our purposes
(such adjustment is needed since the second order condition  (\ref{cond4}) and the corresponding condition  in \cite{Paul16} are slightly different),
we get
that for possibly different function $A_0(t)$ with $A_0(t) \sim A(t)$, as $ t \to \infty$, and for each $\epsilon>0$, $\delta>0$ there
exists a $t_0$ such that for $t>t_0$, $x\ge 1$,
\begin{equation}\label{pthm2-02}
\left|g_{r,0}\left( \frac{U(tx)}{U(t)}\right)-  x^{\gamma r} - r x^{\gamma  r} A_0(t) f_{\rho}(x) \right| \le \epsilon r x^{\gamma r +\rho+\delta} \left|A_0(t)\right|,
\end{equation}
where $f_{\rho}(x)$ is defined in  (\ref{f_r}).
It is well-known that a similar Potter's type bounds hold for the logarithmic function, see, e.g., inequalities (3.2.7) in \cite{Haanbook}. Namely, for a possibly different function $A_1(t)$ with $A_1(t) \sim A(t)$, as $ t \to \infty$, and for each $\epsilon>0$, $\delta>0$ there
exists  $t_1$ such that for $t>t_1$, $x\ge 1$,
\begin{equation}\label{pthm2-03}
\left|\ln\left( \frac{U(tx)}{U(t)}\right)- \gamma \ln(x)-  A_1(t) f_{\rho}(x) \right| \le \epsilon  x^{\rho+\delta} \left|A_1(t)\right|.
\end{equation}
Let $\tilde{t}=\max\{t_0, t_1\}$. By multiplying inequalities (\ref{pthm2-02}) and (\ref{pthm2-03}) we get that for $t>\tilde{t}$, $x\ge 1$,
\begin{eqnarray*}
  &&\bigg| g_{r,1} \left( \frac{U(tx)}{U(t)}\right) - g_{r,0}\left( \frac{U(tx)}{U(t)}\right) \left\{ \gamma \ln(x)+  A_1(t) f_{\rho}(x) \right\}
  \\
  && \quad -  \ln\left( \frac{U(tx)}{U(t)}\right) x^{\gamma  r}\left\{ 1 + r  A_0(t) f_{\rho}(x)\right\}
 \\
&& \quad +
x^{\gamma  r}\left\{ 1 + r  A_0(t) f_{\rho}(x)\right\} \left\{ \gamma \ln(x)+  A_1(t) f_{\rho}(x) \right\} \bigg| \\
&& \ \le \epsilon^2 r x^{\gamma r +2\rho+2\delta} \left|A_0(t)\right| \left|A_1(t)\right|.
\end{eqnarray*}
Suppose that $\tilde{t}$ is large enough that for $t>\tilde{t}$, $x\ge 1$,
$$
\gamma \ln(x)+  A_1(t) f_{\rho}(x) >0, \quad 1 + r  A_0(t) f_{\rho}(x) >0.
$$
Then, by applying inequalities (\ref{pthm2-02}) and (\ref{pthm2-03}) one more time we obtain
\begin{equation}\label{pthm2-04}
-\epsilon d_1(x,t) \le  g_{r,1} \left( \frac{U(tx)}{U(t)}\right)
-b_1(x) -c_1(x,t)
\le \epsilon  d_1(x,t),
\end{equation}
where
\begin{eqnarray*}
b_1(x) &=& \gamma x^{\gamma r} \ln(x), \\
c_1(x,t) &=& \left(A_1(t) +\gamma r A_0(t) \ln(x)\right)x^{\gamma r} f_{\rho}(x)+rA_0(t)A_1(t) x^{\gamma r} f_{\rho}^2(x),  \\
d_1(x,t) &=& x^{\gamma r +\rho+\delta} \bigg(\left|A_1(t)\right| \left\{ 1 + r  A_0(t) f_{\rho}(x)\right\} \\
&& \ +
r \left|A_0(t)\right| \left\{ \gamma \ln(x)+  A_1(t) f_{\rho}(x) \right\}+ \epsilon r x^{\rho+\delta} \left|A_0(t)\right| \left|A_1(t)\right|\bigg),
\end{eqnarray*}

To prove two-dimensional Central Limit Theorem (\ref{res2}) we shall use the well-known Cramer-Wald method. Let $\left(\theta_0,\theta_1\right) \in \mathbf{R}^2$. From (\ref{pthm2-02}) and (\ref{pthm2-04}) we get
\begin{equation}\label{pthm2-05}
-\epsilon  d(x,t) \le  \theta_0 g_{r,0} \left( \frac{U(tx)}{U(t)}\right)+\theta_1 g_{r,1} \left( \frac{U(tx)}{U(t)}\right)
-b(x) -c (x,t)
\le \epsilon  d(x,t),
\end{equation}
where
\begin{eqnarray*}
b(x) &=& \theta_0 x^{\gamma r} +\theta_1 b_1(x), \\
c(x,t) &=& \theta_0 r A_0(t) x^{\gamma r} f_{\rho}(x) +\theta_1 c_1(x,t),  \\
d(x,t) &=& \left| \theta_0\right| r x^{\gamma r+\rho+\delta} |A_0(t)| + \left| \theta_1\right|d_1(x,t).
\end{eqnarray*}
We claim that
\begin{equation}\label{pthm2-06}
\frac{1}{\sqrt{k}} \sum_{i=0}^{k-1} d\left( \frac{Y_{n-i,n}}{Y_{n-k,n}}, Y_{n-k,n}\right) {\buildrel \rm P \over \to}
\frac{r \left| \theta_0 \mu\right|}{d_r(1)-\rho-\delta}+
\frac{(1-\rho-\delta)\left| \theta_1 \mu\right|}{(d_r(1)-\rho-\delta)^2}
, \quad n \to \infty.
\end{equation}
From Lemma 1 in \cite{Paul16} we know that if $\nu<1$, $\nu \not =0$, then
\begin{equation}\label{pthm2-07}
\frac{1}{k} \sum_{i=0}^{k-1} g_{\nu,0}\left(\frac{Y_{n-i,n}}{Y_{n-k,n}}\right) {\buildrel \rm P \over \to} \frac{1}{1-\nu}, \quad n \to \infty.
\end{equation}
Similarly one can prove
\begin{equation}\label{pthm2-08}
\frac{1}{k} \sum_{i=0}^{k-1} g_{\nu,1}\left(\frac{Y_{n-i,n}}{Y_{n-k,n}}\right) {\buildrel \rm P \over \to} \frac{1}{(1-\nu)^2}, \quad n \to \infty.
\end{equation}
The relation
\begin{equation} \label{pthm2-09}
\sqrt{k} A \left(Y_{n-k,n}\right)  {\buildrel \rm P \over \to} \mu, \quad n \to \infty,
\end{equation}
where $\mu$ is the same as in (\ref{cond5}), is proved in \cite{Paul16}. Now, by combining
(\ref{pthm2-07})-(\ref{pthm2-09}) one can obtain (\ref{pthm2-06}).

Taking into account (\ref{pthm2-06}), substituting the values of $t$ and $x$ from  (\ref{pr01-0}) into (\ref{pthm2-05}) and performing summation
over $i=0,1,\dots, k-1$  we get distributional representation
\begin{eqnarray}
 &&\sqrt{k} \left\{\theta_0 \left(G_n(k,r,0)-\xi(r,0)\right)+\theta_1\left( G_n(k,r,1)-\xi(r,1) \right)\right\} \nonumber \\
\label{pthm2-10}
 && \quad {\buildrel \rm d \over =} \sqrt{k} B_n(k,r)+\sqrt{k} C_n(k,r)+o_p(1),
\end{eqnarray}
where
\begin{eqnarray*}
   B_n(k,r) &=& \frac{1}{k} \sum_{i=0}^{k-1} \left\{ b\left( \frac{Y_{n-i,n}}{Y_{n-k,n}}\right) - \theta_0\xi(r,0)-\theta_1\xi(r,1) \right\}, \\
   C_n(k,r) &=& \frac{1}{k} \sum_{i=0}^{k-1} c\left(\frac{Y_{n-i,n}}{Y_{n-k,n}}, Y_{n-k,n}\right).
\end{eqnarray*}
By applying relations (\ref{pthm2-07})-(\ref{pthm2-09}) one more time, we find
\begin{equation} \label{pthm2-11}
\sqrt{k} C_n(k,r)  {\buildrel \rm P \over \to} \mu\left(\theta_0 \nu^{(1)}(r)+ \theta_1 \nu^{(2)}(r)\right),
\quad n \to \infty,
\end{equation}
where $\nu^{(j)}(r)$, $j=1,2$ are defined in (\ref{nu}).
By using well-known  R\'{e}nyi's representation (see e.g. Section 2 in \cite{Paul16} for details) we obtain
\begin{equation} \label{pthm2-12}
\sqrt{k} B_n(k,r)  {\buildrel \rm d \over =} \tilde{B}_n(k,r),
\end{equation}
where
$$
\tilde{B}_n(k,r)= \frac{1}{\sqrt{k}} \sum_{i=0}^{k-1} \left\{ \theta_0\left( g_{\gamma r,0}(Y_i) - \xi(r,0)\right)+
\theta_1 \gamma \left( g_{\gamma r,1}(Y_i) - \frac{1}{d_r^2(1)}\right) \right\}.
$$
Keeping in mind equality (\ref{id01})
one can deduce that
quantity $\tilde{B}_n(k,r)$ presents normalized sum of zero mean random variables, which are i.i.d.. Moreover, under assumption $\gamma r<1/2$,
\begin{eqnarray*}
\E \left\{ \theta_0\left( g_{\gamma r,0}(Y_i) - \xi(r,0)\right)+
\theta_1 \gamma \left( g_{\gamma r,1}(Y_i) - \frac{1}{d_r^2(1)}\right) \right\}^2  = \theta_0^2 s_1^2(r) +2 \theta_0 \theta_1 s_{12}(r)+\theta_1^2 s_1^2 (r),
\end{eqnarray*}
where $s_1^2(r)$, $s_2^2(r)$ and $s_{12}(r)$ are defined in (\ref{s}). Thus, applying Lindeberg-L\'{e}vy central limit theorem we get the relation
$\tilde{B}_n(k,r) {\buildrel \rm d \over \to} \theta_0 W^{(1)}+\theta_1 W^{(2)}$. This, together with (\ref{pthm2-12}) gives
\begin{equation} \label{pthm2-13}
\sqrt{k} B_n(k,r)  {\buildrel \rm d \over \to} \theta_0 W^{(1)}+\theta_1 W^{(2)}, \quad n \to \infty.
\end{equation}
Applying  Theorem 3.9 in \cite{Billingsley}, from  (\ref{pthm2-11}) and (\ref{pthm2-13}) we get
$$
\sqrt{k} \left(B_n(k,r),  C_n(k,r)\right) {\buildrel \rm d \over \to} \left(\theta_0 W^{(1)}+\theta_1 W^{(2)}, \mu\left(\theta_0 \nu^{(1)}(r)+ \theta_1 \nu^{(2)}(r)\right) \right), \quad n \to \infty.
$$
 Continuous Mapping Theorem gives us the relation
$$
\sqrt{k} \left(B_n(k,r)+C_n(k,r)\right) {\buildrel \rm d \over \to}
\theta_0 \left(W^{(1)}+ \mu \nu^{(1)}(r)\right)+ \theta_1 \left(W^{(2)}+ \mu \nu^{(2)}(r)\right), \quad n \to \infty.
$$
The last relation together with (\ref{pthm2-10}) gives (\ref{res2}). Theorem \ref{thm2} is proved. \halmos

\bigskip

\noindent {\it Proof of the Theorem \ref{thm3}.} In the case $j=1$ the proof of the relation (\ref{res3}) can be found in \cite{Paul16} (proof of the Corollary 1) or
in \cite{Beran} (proof of Theorem 2). But the  asymptotic normality of all estimators $ \hat{\gamma}^{(j)}_n(k,r), \ j=1,2,3$
can be obtained in a unified way, expressing these estimators as functions of statistics $G_n(k,r,0)$ and $G_n(k,r,1)$, and then  combining Theorem \ref{thm1}, Theorem \ref{thm2} and Continuous mapping Theorem. Namely, it is easy to see that
$$
\hat{\gamma}^{(1)}_n(k,r)-\gamma= \frac{(1-\gamma r)\left(G_n(k,r,0)-\xi_{r,0}\right)}{r G_n(k,r,0)}.
$$
For the estimator $\hat{\gamma}^{(2)}_n(k,r)$ we have
\begin{eqnarray*}
 \hat{\gamma}^{(2)}_n(k,r)-\gamma
 = \frac{2(1-\gamma r)G_n(k,r,1)-\gamma - \gamma \sqrt{4r G_n(k,r,1)+1} }{2r G_n(k,r,1)+1+\sqrt{4r G_n(k,r,1)+1}}.
\end{eqnarray*}
Multiplying  the numerator and denominator of the right hand side by $2(1-\gamma r)G_n(k,r,1)-\gamma + \gamma \sqrt{4r G_n(k,r,1)+1}$, we get
\begin{eqnarray*}
 && \hat{\gamma}^{(2)}_n(k,r)-\gamma \\
&& \quad = \frac{4(1-\gamma r)^2 \left(G_n(k,r,1)-\xi_{r,1}\right)\left(G_n(k,r,1)+\xi_{r,1}\right)-4\gamma \left(G_n(k,r,1)-\xi_{r,1}\right) }{2r G_n(k,r,1)+1+\sqrt{4r G_n(k,r,1)+1}} \times \\
&& \quad \ \times \frac{1}{ 2(1-\gamma r)G_n(k,r,1)-\gamma + \gamma \sqrt{4r G_n(k,r,1)+1}}.
\end{eqnarray*}
For the third estimator the following representation holds:
\begin{eqnarray*}
\hat{\gamma}^{(3)}_n(k,r)-\gamma &=& \frac{r(1-\gamma r)\left(G_n(k,r,1)-\xi_{r,1} \right)-\left(G_n(k,r,0)-\xi_{r,0} \right)}{r^2 G_n(k,r,1)}.
\end{eqnarray*}
As it was said, it remains to combine Theorem \ref{thm1}, Theorem \ref{thm2} and Continuous Mapping Theorem, and  we deduce  (\ref{res3}) with $1 \le j \le 3$. For example, we have
$$
\hat{\gamma}^{(3)}_n(k,r)-\gamma =f\left (\sqrt{k} (G_n(k,r,1)-\xi(r,1)),\sqrt{k} (G_n(k,r,0)-\xi(r,0)), G_n(k,r,1)\right )
$$
with
$$
f(x,y,z)=\frac{r(1-\g r)x-y}{r^2 z}.
$$
From Theorems \ref{thm1} and \ref{thm2} and Theorem 3.9 from \cite{Billingsley} we have
$$
\left ( \sqrt{k} \left(G_n(k,r,i)-\xi(r,i)\right ), i=0,1, G_n(k,r,1) \right) {\buildrel \rm d \over \to}
    \left(W^{(i)}+\mu \nu^{(i)}(r),i=1,2, \xi(r,1) \right),
$$
and now we apply Continuous Mapping Theorem (Theorem 2.7 in \cite{Billingsley}).
\halmos

As it was mentioned at the end of Introduction, there is possibility to use general approach in proving asymptotic normality of the introduced estimators, suggested in \cite{Drees}. Let $F_n$ stands for the empirical d.f. based on the sample $X_1, \dots, X_n$ and let the empirical tail quantile function is defined as
$$
Q_n(t):= F_n^{-1}\left (1-\frac{k_n}{n}t\right )=X_{n-[k_nt], n}, \quad  t\in [0,1].
$$
Then almost all known estimators of the extreme value index that are based on some part of largest order statistics can be written as some functional
$T$ (defined on some functional space), applied to $Q_n$. Then, having estimator written  as $T(Q_n)$, the idea in \cite{Drees} is to use invariance principle for the process $Q_n$ in the functional space on which $T$ is defined and then requiring some smoothness (Hadamard differentiability in linear topological space) one can try to derive asymptotic normality of the estimator under consideration. For estimators, considered in the paper, it is possible to use this scheme, but the functionals are quite complicated. For example, one can write $\hat{\gamma}^{(3)}_n(k,r)=T_{GMR}(Q_n)$, but
$$
T_{GMR}(z)=f(T_0(z), T_1(z), T_2(z) ),
$$
$$
   f(x,y,v)=\frac{ry-x+1}{r^2y}, \ {\rm for} \ \ r\ne 0  \quad {\rm and} \ \ f(x,y,v)=\frac{z}{2y}, \ {\rm for} \ \ r= 0,
$$
and
$$
T_i(z)=\int_0^1 g_{r,i}\left (\frac{z(t)}{z(1)} \right )dt, \ \ i=0,1,  \quad T_2(z)=\int_0^1 g_{0,2}\left (\frac{z(t)}{z(1)} \right )dt.
$$
 For this complicated functional we must prove Hadamard differentiability on some linear topological space (to choose the appropriate space  is also non trivial problem,  usual $D[0,1]$ space with Skorokhod topology is not a linear topological  space). Thus, it seems that our approach is much more simple and we do not need more restrictive conditions (such that appears in \cite{Drees}), since instead of invariance principle for tail quantile function we prove two-dimensional CLT for two particular statistics and then apply continuous mapping theorem in $R^3$.

\noindent {\it Proof of the Theorem \ref{thm4}}. The expression of  $B^{(1)}(r)$ is given  in \cite{Beran}, and at first we followed the pattern of the proof in \cite{Beran}, but, as it was noticed in Introduction, there is more simple proof.

From the expression (\ref{robust1}), given in the Introduction, we see that, in order to find $B_n^{(\ell)} (k,r)$,
 it is sufficient to find the limit $\lim_{x \to \infty}\hat{\gamma}^{(\ell)}_n(k,r; X_1, \dots, X_{n-1},x)$, since $\hat{\gamma}^{(\ell)}_{n-1}(k-1,r; X_1, \dots, X_{n-1}) {\buildrel \rm P \over \to} \g$, as $n\to \infty$ and (\ref{cond3}) holds. For all three estimators calculations are simple and similar, therefore we demonstrate the proof for the estimator $\hat{\gamma}^{(2)}_n(k,r,x)$, having  the most complicated expression.
It is clear that, for sufficiently large value of $x$, $\hat{\gamma}^{(\ell)}_n(k,r; X_1, \dots, X_{n-1},x)$ can be written as
\begin{equation}\label{robust}
\frac{2\left (h(x)+b \right )}{2r\left (h(x)+b \right )+1+\sqrt{4r \left (h(x)+b \right )+1} },
\end{equation}
where $h(x)=g_{r,1}\left(x/X_{n-k,n} \right)$ and $b$ is the sum of the rest summands from statistic $G_n(k,r, 1)$ and does not depend on $x$.
If $r>0$, then $h(x)\to \infty$ and the limit of the quantity in (\ref{robust}) is $1/r$, while  for $r<0$ \ $h(x)\to 0$, and the limit in (\ref{robust})  is
$$
\frac{2b }{2rb +1+\sqrt{4r b +1} }
$$
and this expression almost (this word is used for the reason that in the above written expression the sum is divided by  $k$, while for complete coincidence division should be by $k-1$ ) coincides with $\hat{\gamma}^{(\ell)}_{n-1}(k-1,r; X_1, \dots, X_{n-1})$, therefore, passing to the limit as $n\to \infty$ we get in limit $B^{(2)}=0$. In the case $r=0$ only nominator contains function $h(x)=g_{0,1}\left(x/X_{n-k,n} \right)$  which grows unboundedly, therefore, we get infinite value for $B_n^{(\ell)} (k,r)$.
\halmos

%
%


\end{document}